\documentclass[journal]{IEEEtran}
\usepackage{graphicx}
\usepackage{amsfonts,amssymb}
\usepackage{amsmath}
\usepackage{amsthm}
\usepackage{cases}
\usepackage{mathrsfs}
\usepackage{mathabx}
 \usepackage{xcolor}
\usepackage{savesym}
\savesymbol{AND}
\usepackage{algorithm}
\usepackage{algorithmic}

\usepackage[colorlinks,
linkcolor=blue,
anchorcolor=blue,
citecolor=red
]{hyperref}
\usepackage{cite}
\usepackage{subcaption}

\newtheorem{theorem}{Theorem}
\newtheorem{lemma}{Lemma}
\newtheorem{remark}{Remark}
\newtheorem{pro}{Proposition}
\newtheorem{coro}{Corollary}

\newtheorem{definition}{Definition}
\newtheorem{exam}{Example}
\newtheorem{assu}{Assumption}

\newcommand*{\blue}{ }

\newcommand{\rank}{\mathrm{Rank}}

\begin{document}
	
	%\begin{frontmatter}
	\title{Identifiability in Dynamic Acyclic Networks with Partial Excitations and Measurements}

	\author{Xiaodong~Cheng, 
		Shengling Shi,
		Ioannis Lestas, and
		Paul M.J. Van den Hof 
		%<-this % stops a space <-this % stops a space
		\thanks{This work is supported by 
			the European Research Council (ERC), Advanced Research Grant SYSDYNET (No. 694504) and Starting Grant HetScaleNet (No. 679774).}%
		\thanks{
			X. Cheng is with Mathematical and Statistical Methods Group (Biometris), Department of Plant Science, Wageningen University \& Research,
			6700 AA Wageningen, The Netherlands.
			{\tt\small xiaodong.cheng@wur.nl}}
		\thanks{
			S. Shi is  with the Delft Center for Systems and Control, Delft University of Technology, The Netherlands. 
			{\tt\small s.shi-3@tudelft.nl}}
		\thanks{
			I. Lestas is with the Control Group,  Department
			of Engineering, University of Cambridge, CB2 1PZ, United Kingdom.
			{\tt\small icl20@cam.ac.uk}}
		\thanks{P. M. J. Van den Hof are with the Control Systems Group, Department of Electrical Engineering,
			Eindhoven University of Technology, 
			5600 MB Eindhoven,
			The Netherlands.
			{\tt\small P.M.J.vandenhof@tue.nl}}%
		
	}
	
	\maketitle

	%\begin{IEEEkeyword}                           
	% 	System identification, Network Identifiability, Sensor allocation, Dynamic networks,   Acyclic Networks     
	%\end{IEEEkeyword} 
	
	\begin{abstract}
		
		This paper deals with dynamic networks in which the causality relations between the vertex signals are represented by linear time-invariant transfer functions (modules).
		Considering an acyclic network where only a subset of its vertices are measured and a subset of the vertices are excited, we explore conditions under which all the modules are identifiable on the basis of measurement data. Two sufficient conditions are presented in the paper, 
		where the first condition concerns an identifiability analysis that needs to be performed for each vertex, while the second condition, based on the concept of \textit{tree/anti-tree covering}, results from a graphical synthesis approach to allocate actuators and sensors in an acyclic network for achieving generic identifiability.
		%  
		%  where the former one can be applied to perform an identifiability analysis for a full network, while the latter one, based on the concept of \textit{tree/anti-tree covering}, induces a synthesis approach to allocate actuators and sensors in an acyclic network for generic identifiablity.
	\end{abstract}

	%\end{frontmatter}
	\section{Introduction} 
	
	Dynamic networks have appeared in a wide range of technological applications including biochemical reaction systems, decentralized manufacturing processes, and smart power grids. In recent years, considerable attention from the systems and control domain has been devoted to address data-driven modeling problems in dynamic networked systems with large-scale interconnection structure, see some of the representative works in \cite{gonccalves2008LTInetworks,pare2013necessity,chetty2017necessary,materassi2012problem,nabi2012network,adebayo2012dynamical,hayden2016sparse,bazanella2017identifiability,hendrickx2018identifiability,Yu2018identification,paul2013netid,weerts2018identifiability,henk2018identifiability}. From an identification perspective, the modeling problems in a network setting can be viewed as an extension of classical system identification on the basis of open-loop data or closed-loop data, towards an identification problem with data collected from distributed dynamical systems operating under a network interconnection. This leads to a framework where signals are regarded as vertices in a network, and the causal dependencies among the signals, typically represented by proper transfer functions (modules), are considered
	as directed edges.
	
	Different data-driven modeling
	problems can be formulated in such a framework. For example, how to detect the topology of a network from measurement data  \cite{materassi2010topology,chiuso2012bayesian,materassi2012problem,yeung2015global,shi2019bayesian}?
	How to consistently estimate selected local dynamics within a network \cite{paul2013netid,dankers2016identification,gevers2018practical,materassi2020signal,ramaswamy2021local,weerts2020abstractions}? What are the conditions required to identify the dynamics of an entire network or a subnetwork \cite{weerts2015sysid,bazanella2017identifiability,weerts2018identifiability,hendrickx2018identifiability,henk2018identifiability,cheng2019allocation,shi2020subnetworks}? In this paper, we focus on the last problem, where  \textit{identifiability} of a full network based on a subset of external excitation signals and measured vertex signals is of particular interest. This concept plays a key role in data-driven modeling of dynamic networks, as it essentially reflects if we are able to acquire a unique network model on the basis of measurement data. In the literature, there are two classes of network identifiability, namely, \textit{global identifiability} \cite{henk2018identifiability,weerts2018identifiability} that requires \textit{all} the models in the set to be distinguishable, and \textit{generic identifiability} \cite{bazanella2017identifiability,hendrickx2018identifiability,weerts2018single,Bazanella2019CDC}, which means that \textit{almost all} models in the considered model set can be distinguished. Recently, the concept of \textit{local generic identifiability} has been introduced in \cite{legat2020local,legat2021path}, which can be viewed as a notion that is weaker than generic identifiability.

	In this paper, we will focus on global and generic identifiability notions, which can be addressed in two problem settings, namely, dynamic networks with full excitation/measurement, and the ones with partial excitation/measurement.
	%In e.g., \cite{hendrickx2018identifiability,henk2018identifiability}, all vertices are excited by external signals, while only a subset of vertices is measured, and in e.g., \cite{weerts2018identifiability,weerts2018single}, all internal variables are supposed to be measured, while only partial vertices are excited or influenced by noises. 
	A number of conditions have been derived for (generic) identifiability of a full dynamic network, see e.g., \cite{bazanella2017identifiability,hendrickx2018identifiability,henk2018identifiability,weerts2018identifiability,cheng2019allocation}, which are based on network topology i.e., how the vertices are interconnected. While the identifiability analysis in \cite{bazanella2017identifiability,hendrickx2018identifiability,henk2018identifiability} is performed under the assumption that all the vertices are excited by sufficiently rich external signals, these works lead to attractive path-based conditions for checking network identifiability. Depending on whether \textit{generic identifiability} or \textit{global identifiability} is considered, these conditions are interpreted in terms of the existence of \textit{vertex disjoint paths}\cite{bazanella2017identifiability,hendrickx2018identifiability} or \textit{constrained vertex disjoint paths}\cite{henk2018identifiability} from the out-neighbors of each vertex to the measured vertices in a network.
	In contrast, the result in \cite{weerts2015sysid,weerts2018identifiability} is developed for the setting where only a subset of vertices are affected by external signals including noises and excitation signals manipulated by users, while all the vertices are measured. Moreover, network identifiability is defined as a property of a parameterized model set, instead of a property of a single network as done in \cite{bazanella2017identifiability,hendrickx2018identifiability,henk2018identifiability}. 
	%Such a definition allows for a more flexible choice of model sets that can include nonparameterized network modules and known correlations between disturbance signals in a network.
	In \cite{weerts2015sysid,weerts2018identifiability}, global network identifiability is characterized by the rank property of a certain transfer matrix that is determined by the presence and location of external signals,
	the correlation structure of disturbances, and the topology of parametrized modules. This rank condition has been further studied in \cite{weerts2018single} leading to generic identifiability of a network model set, which is equivalent to a vertex-disjoint path condition dual to the ones in \cite{bazanella2017identifiability,hendrickx2018identifiability}. With the same network setting as in \cite{weerts2015sysid,weerts2018identifiability,weerts2018single}, the work in \cite{cheng2019CDC,cheng2019allocation} provides a new characterization for generic identifiability using the concept of \textit{disjoint pseudo-tree covering}, from which a
	graphical tool has been developed for allocating external excitation signals so as to warrant generic identifiability.

	While different measurement and	excitation schemes are considered, most works in the existing literature require all the vertices to be either simultaneously excited or simultaneously measured. The work in \cite{Bazanella2019CDC} relaxes this requirement and investigates the identifiability condition with partial excitation and partial measurement, although these conditions hold in some special cases. For example, {\blue some graph-based results are provided for networks with special tree and cycle topologies, and some require that all the out-neighbors of a vertex are excited, or all the in-neighbors of a vertex are measured.} In contrast, \cite{shi2020partial}  provides a sufficient condition for generic identifiability of a single module in a network, and the condition does not have specific requirements as in  \cite{Bazanella2019CDC} and is applicable to the more general case.

	In line with the network setting of \cite{Bazanella2019CDC,shi2020partial}, the current paper aims to derive general conditions for identifiability in full dynamic networks, where only partial excitation and partial measurement signals are available. 
	First, we introduce the concept of transpose networks  
	%in which excited (measured) vertices in the original network become measured (excited), and the   edges directions are reversed compared to the original network. 
	and show that identifiability of the model sets of a transpose network and that of its corresponding original network are equivalent, which builds a connection between the two settings in  \cite{bazanella2017identifiability,hendrickx2018identifiability,henk2018identifiability} and \cite{weerts2018identifiability,weerts2018single,cheng2019allocation}. Although similar duality discussions have been seen in \cite{hendrickx2018identifiability,cheng2019allocation}, the analysis therein is performed in the setting that the vertices in a network are fully excited or measured. 
	In this paper, we further generalize the duality analysis to dynamic networks with partial excitation and measurement by using transpose networks.
	Second, we provide new insight into the identifiability analysis of full networks with  \textit{acyclic topology} in the most general setting of \textit{partial excitation and partial measurement}. {\blue Compared to \cite{mapurunga2022excitation}, which presents necessary conditions for identifiability of  acyclic networks, 
		as our main contribution,
		we present sufficient conditions to determine both global and generic identifiability of acyclic networks.}	The first condition takes advantage of the hierarchical structure in acyclic networks to formulate rank conditions for every vertex, which can also be reinterpreted as path-based conditions and hence lead to a vertex-wise check for identifiability of a given acyclic network. {\blue Furthermore}, 
	% to further solve the problem of actuator and sensor allocation, 
	we present an alternative condition for checking generic identifiability of model sets of acyclic networks on the basis of a graphical concept, called \textit{disjoint tree/anti-tree covering}. This condition then leads to an algorithmic procedure that achieves generic identifiability of the model set of an acyclic network requiring as few as possible excited and measured vertices. 
	
	The rest of this paper is organized as follows. In Section~\ref{sec:preliminaries}, we recap some basic notations in graph theory and introduce the dynamic network model. Section~\ref{sec:TransposeNet} presents the results on transpose networks, and Section~\ref{sec:Acyclic} provides the main results of the paper on identifiability in acyclic networks, and a graph-theoretic approach is also presented for the allocation of actuators and sensors. Finally, concluding remarks are made in Section~\ref{sec:conclusion}. The proofs of the technical results are presented in the appendix.

	\section{Preliminaries and Dynamic Network Setup}
	
	\label{sec:preliminaries}
	
	\subsection{Notation}
	Denote $\mathbb{R}$ as the set of real numbers and $\mathbb{R}(q)$ as the rational function field over $\mathbb{R}$ with the variable $q$. $I_n$ is the identity matrix of dimension $n$. The cardinality of a set $\mathcal{V}$ is represented by $\lvert \mathcal{V} \rvert$. $A_{ij}$ denotes the $(i,j)$ entry of a matrix $A$, and more generally, $[A]_{\mathcal{U},\mathcal{V}}$ denotes the submatrix of $A$ that consists of the rows and columns of $A$ indexed by two positive integer sets $\mathcal{U}$ and $\mathcal{V}$, respectively. Furthermore, $[A]_{\mathcal{U}, \star}$ and $[A]_{\star, \mathcal{V}}$ represent the matrix containing the rows of $A$  indexed by the set $\mathcal{U}$ and the matrix containing the columns of $A$  indexed by the set $\mathcal{V}$, respectively.

	The normal rank of a transfer matrix $A(z)$ is denoted by $\rank(A(z))$, and $\rank(A(z)) = r$ if the rank of $A(z)$ is equal to $r$ for almost all values of $z$. Furthermore, let $A(z,\theta)$ be a parameterized transfer matrix with the parameters $\theta \in \Theta$, then the generic rank of $A(z,\theta)$, denoted by $\mathrm{gRank}(A)$, is the maximum normal rank of $A(z,\theta)$ for all $\theta \in \Theta$ \cite{vanderWoude1991graph}.

	\subsection{Graph Theory}
	
	The topology of a dynamic network is characterized by a graph $\mathcal{G}$ that consists of a finite and nonempty vertex set  $\mathcal{V}: = \{1, 2, ... , L\}$ and an edge set $\mathcal{E} \subseteq \mathcal{V} \times \mathcal{V}$.
	In a directed graph each element in $\mathcal{E}$ is an ordered pair of elements of $\mathcal{V}$. If $(i,j) \in \mathcal{E}$, we say that the edge is incident from vertex $i$ to vertex $j$, the vertex $i$ is the \textit{in-neighbor} of $j$, and $j$ is the \textit{out-neighbor} of $i$. Let $\mathcal{N}_j^-$ and $\mathcal{N}_j^+$ be the sets that collect all the in-neighbors and out-neighbors of vertex $j$, respectively.
	
	A graph $\mathcal{G}$ is called \textit{simple}, if $\mathcal{G}$ does not contain self-loops (i.e., $\mathcal{E}$ does not contain any edge of the form $(i,i)$, $\forall~i\in \mathcal{V}$), and there exists only one directed edge from one vertex to each out-neighbor. In a simple graph, a directed \textit{path} connecting vertices
	$i_0$ and $i_n$ is a sequence of edges of the form $(i_{k-1}, i_k)$, $k = 1, \cdots, n$, and every vertex appears at most once on the path. If there is a directed path from vertex $i$ to $j$, we say that $j$ is \textit{reachable} from $i$, {\blue and $j$ is a \textit{descendant} of $i$}. Two directed paths are \textit{vertex-disjoint} if they do not share any common vertex, including the start and the end vertices. 
	In a simple directed graph $\mathcal{G}$, we denote $b_{\mathcal{U} \rightarrow \mathcal{Y}}$ as the maximum number of mutually vertex-disjoint
	paths from $\mathcal{U} \subseteq \mathcal{V}$ to $\mathcal{Y} \subseteq \mathcal{V}$. 
	%Furthermore, a set of $m$  vertex-disjoint paths from $\mathcal{U}$ to $\mathcal{Y}$ is \textit{constrained} if it  is unique, and we denote $\hat{b}_{\mathcal{U} \rightarrow \mathcal{Y}}$ as the maximum number of constrained vertex-disjoint paths from $\mathcal{U}$ to $\mathcal{Y}$.
	%A directed simple graph $\mathcal{G}$ is \textit{connected} if the underlying undirected graph $\mathcal{G}_\mathrm{u}$ obtained by replacing all directed edges of $\mathcal{G}$ with undirected edges is connected, i.e., in $\mathcal{G}_\mathrm{u}$, there is an undirected path between any pair of vertices.

	Let $\mathcal{U}$ and $\mathcal{Y}$ be two vertex subsets in a directed graph $\mathcal{G}$. A set of vertices $\mathcal{D}$ in $\mathcal{G}$ is a \textit{disconnecting set} from $\mathcal{U}$ to $\mathcal{Y}$, if all the directed paths from $\mathcal{U}$ to $\mathcal{Y}$ in $\mathcal{G}$  pass through $\mathcal{D}$. Roughly speaking, removing all the vertices in a disconnecting set $\mathcal{D}$ from $\mathcal{U}$ to $\mathcal{Y}$ will lead to the absence of directed paths from $\mathcal{U}$ to $\mathcal{Y}$. Note that it is allowed for $\mathcal{D}$ to include vertices in the sets $\mathcal{U}$ and $\mathcal{Y}$.

	In this paper, we concentrate on simple \textit{acyclic} graphs, meaning that a vertex does not reach itself via a directed path.
	In a simple acyclic graph $\mathcal{G}$, a \textit{source} is a vertex without any in-neighbors, and likewise, a \textit{sink} is a vertex without any out-neighbors. 
	
	%We denote ${\mathcal S_\mathrm{ou}}(\mathcal{G})$ and ${\mathcal{S}_\mathrm{in}}(\mathcal{G})$ as the sets of the sources and sinks of $\mathcal{G}$, respectively.

	\subsection{Dynamic Network Model}
	
	Consider a dynamic network whose topology is captured by a simple directed graph $\mathcal{G} = (\mathcal{V}, \mathcal{E})$ with a vertex set
	$\mathcal V = \{1, 2, ... , L\}$ and edge set $\mathcal E \subseteq \mathcal V \times \mathcal V$. Following the basic setup in  \cite{paul2013netid,Bazanella2019CDC}, each vertex describes an internal variable $w_j(t) \in \mathbb{R}$, and 
	a compact form of the overall network dynamics is 
	\begin{align} \label{eq:net}
		w(t) & = G(q) w(t) +  R r(t) + v(t), \\
		y(t) & = C w(t), \nonumber
	\end{align}
	where $G(q)$ is a rational matrix in the delay operator $q^{-1}$, with zero diagonal elements, $v(t) \in \mathbb{R}^L$ is a vector of zero-mean, stationary stochastic process noises, and $R \in \mathbb{R}^{L\times K}$ and $C \in \mathbb{R}^{N\times L}$ are binary selection matrices. All the internal variables are stacked in the vector $w(t): = [
	w_1(t) \ w_2(t)\ \cdots \ w_L(t)
	]^\top$. $G(q)$ is the system matrix which is a matrix of transfer functions with zero diagonal elements. The $(i,j)$-th entry of $G(q)$, denoted by $G_{ij}(q) \in \mathbb{R}(q)$, indicates the transfer operator from vertex $j$ to vertex $i$, and it is represented by an edge $(j,i) \in \mathcal{E}$ in graph $\mathcal{G}$.
	
	Let $\mathcal{R} \subseteq \mathcal{V}$ and $\mathcal{C} \subseteq \mathcal{V}$ be the vertices that are excited and measured, respectively, and $K = |\mathcal{R}|$ and $N = |\mathcal{C}|$. The signals $r(t) \in \mathbb{R}^K$ and $y(t) \in \mathbb{R}^N$ are the external excitation and measurement signals with $R = [I_L]_{\star, \mathcal{R}}$, $C = [I_L]_{\mathcal{C},\star}$ binary matrices indicating which vertices are excited and measured, respectively. Particularly, if a vertex $k$ is unexcited, i.e. $k \in \mathcal{V} \setminus \mathcal{R}$, then the $k$-th row of $R$ is zero, and analogously, if a vertex is unmeasured, then the corresponding column of $C$ is zero.
	
	\begin{assu} \label{Assum}
		Throughout the paper, we consider a dynamic network \eqref{eq:net} with the following properties.
		\begin{enumerate}
			\item The network \eqref{eq:net} is \textit{well-posed} and stable, i.e. $(I - G(q))^{-1}$ is proper
			and stable;
			
			\item All the entries of $G(q)$ are proper transfer operators.
			
		\end{enumerate}
	\end{assu}
	
	The above assumptions are standard in the identification of dynamic networks, see e.g., \cite{paul2013netid,hendrickx2018identifiability,henk2018identifiability,Bazanella2019CDC}, and ensure the properness and stability of the mapping from $r$ to $w$.
	
	For specifying the notion of identifiability, we first need to introduce a network model set. 
	%Based on the model setting in \eqref{eq:net}, a problem of interest concerns if all the dynamics of a network, i.e., parameterized transfer functions in $G(q)$, can be consistently identified from the external excitation signals $r$ and the measured output $y$. 
	Let $M = (G, R, C)$ be a network model of \eqref{eq:net} and 
	\begin{equation} \label{eq:modelset}
		\mathcal{M}: = \{M(q,\theta) = (G(q,\theta), R, C), \theta \in \Theta\}
	\end{equation} 
	be a network model set with  parameterized models $M(q,\theta)$, in which all the nonzero transfer functions in $G(q,\theta)$ are parameterized  independently. The network model set $\mathcal{M}$ inherently contains 
	%{\blue includes}  
	information on the topology of a dynamic network, as well as on 
	%non-parameterized modules,
	the presence and locations of actuators and sensors.
	Then identifiability of the network model set is defined as follows \cite{bazanella2017identifiability,hendrickx2018identifiability,shi2020subnetworks,shi2020partial}.
	\begin{definition}
		\label{defn:netid} 
		Denote the transfer matrix
		\begin{equation} \label{eq:T}
			T(q, \theta) : = (I - G(q,\theta))^{-1}.
		\end{equation} 
		Given the network model set $\mathcal{M}$ in \eqref{eq:modelset}. Consider $\theta_0 \in \Theta$ and the following implication
		\begin{align}\label{eq:implication}
			C T(q, \theta_1) R =  C T(q, \theta_0) R
			\Rightarrow
			G_{ji}(q,\theta_1) = G_{ji}(q,\theta_0) 
		\end{align}
		for all parameter $\theta_1 \in \Theta$.
		Then the module $G_{ji}$ is
		\begin{itemize}
			\item \textbf{identifiable} in $\mathcal{M}$ from the submatrix $[T]_{\mathcal{C},\mathcal{R}}$ 
			if the implication \eqref{eq:implication} holds for all $\theta_0 \in \Theta$;
			\item  
			\textbf{\blue generically identifiable} in $\mathcal{M}$ from the submatrix $[T]_{\mathcal{C},\mathcal{R}}$ if the implication \eqref{eq:implication} holds for almost all $\theta_0 \in \Theta$.
		\end{itemize} 
		The model set $\mathcal{M}$ is ({\blue generically}) identifiable if $G_{ji}$ is ({\blue generically}) identifiable for all $(i,j)  \in \mathcal{E}$.
	\end{definition}
	
	In \cite{henk2018identifiability}, identifiability in Definition~\ref{defn:netid} is also referred to as ``global identifiability''. In this paper, we simply use the term ``identifiability'' for ``global identifiability'' for ease of notation. The term ``almost all'' means the exclusion of parameters that are in a subset of $\Theta$ with Lebesgue measure zero. We refer to  \cite{bazanella2017identifiability,hendrickx2018identifiability} for more details about the notion of generic identifiability. 
	\begin{remark}
		%		Note that the varia$[T]_{\mathcal{C},\mathcal{R}}$, el set $\mathcal{M}$, while the properties of the mapping from $\theta$ to network models is not the focus of this study. 
		Although we consider only $[T]_{\mathcal{C},\mathcal{R}}$, i.e. the transfer from $r$ to $y$, as a basis for  
		%excitation input $r$ in the definition of 
		identifiability, the disturbances $v$ in \eqref{eq:net} can also be taken into account as is done in \cite{cheng2019allocation}. Under some mild assumptions, disturbance inputs play a similar role as excitation inputs,  and therefore the results in this paper can directly be generalized to include $v$ by using the notion of extended graphs (\cite{cheng2019allocation}). For the sake of simplicity, we will not address that generalization in this paper.
	\end{remark}
	
	%Note that the variable $\theta \in \Theta$ is only used for formalizing the model set $\mathcal{M}$, while the properties of the mapping from $\theta$ to network models is not the focus of this study. 
	Identifiability of a  dynamic network model set reflects the ability to distinguish between models in the set $\mathcal{M}$ from measurement data, or more precisely, from the transfer matrix $[T]_{\mathcal{C},\mathcal{R}}$ as described in Definition~\ref{defn:netid}. In this sense, network identifiability essentially depends on the presence and location of external excitation signals $r$ and the selection of measured vertex signals $y$.

	\section{Dynamic Transpose Networks}
	\label{sec:TransposeNet}
	
	%We make the results in this section self-contained and independent of the identifiability analysis in acyclic networks.
	In network identifiability studies, typically two situations are distinguished where either all vertices are measured, or all vertices are excited, while it has been shown that both situations lead to results that are dual, see \cite{hendrickx2018identifiability,cheng2019allocation}. 
	%	In the context of network identifiability, dual settings of dynamic networks have been discussed in \cite{hendrickx2018identifiability,cheng2019allocation}, in which all vertices are either measured or excited. 
	In this section, we present a generalization of the duality results for network identifiability with partially excited and partially measured vertices, which is based on the notion of transpose graphs. This generalization holds for general networks, including both cyclic and acyclic networks. The result in this section will then be applied to derive the main results in Section~\ref{sec:Acyclic}.
	
	Let $\mathcal{G} = (\mathcal{V}, \mathcal{E})$ be a simple directed graph. The \textit{transpose graph} of $\mathcal{G}$ is defined by $\mathcal{G}^\prime: = (\mathcal{V}, \mathcal{E}^\prime)$, where $\mathcal{E}^\prime$ is obtained by reversing the direction of each edge in $\mathcal{E}$. Thus, the term `transpose' is because the adjacency matrix of $\mathcal{G}^\prime$ is the transpose of that of the original graph  $\mathcal{G}$. A dynamic transpose network model can be defined on the basis of  $\mathcal{G}^\prime$  as follows.
	\begin{definition}[Dynamic Transpose Network]
		\label{defn:transpose}
		Consider the dynamical model $M = (G, R, C)$ in \eqref{eq:net} associated with a directed graph $\mathcal{G}$. A transpose network model of $M$ is defined based on the transpose graph $\mathcal{G}^\prime$ as
		\begin{align} \label{eq:nettrans}
			\tilde{M} = (\tilde{G}, \tilde{R}, \tilde{C}),
		\end{align}
		where $\tilde{G} = G^\top$, $\tilde{R} = C^\top$, and $\tilde{C} = R^\top$.
	\end{definition}
	
	The transpose network in \eqref{eq:nettrans} is defined as the conceptually dual dynamical model of \eqref{eq:net}, and it has the same set of vertices with every edge reversed compared to the orientation of the corresponding edge in the original network \eqref{eq:net}. Furthermore, all the excited (measured) vertices in the original network become measured (excited) ones in the transpose network. The relation between the identifiability properties in the original and transpose networks is presented in the following lemma, see its proof in Appendix~\ref{ap:lem:dualty}.
	
	\begin{lemma}[Duality in identifiability] \label{lem:dualty}
		Let $\mathcal{M}^\prime : =  \{\tilde{M}(q,\theta) = (G(q,\theta)^\top, C^\top, R^\top), \theta \in \Theta\}$ be the parameterized model set of the transpose network \eqref{eq:nettrans}. Then, the following statements hold.
		\begin{enumerate}
			\item A module ${G}_{ij}$ is identifiable in $\mathcal{M}$ if and only if the module $\tilde{G}_{ji}$ is identifiable in $\mathcal{M}^\prime$; 
			
			\item The network model set $\mathcal{M}$ is identifiable if and only if the network model set $\mathcal{M}^\prime$ is identifiable.
		\end{enumerate}
		
	\end{lemma}
	
	Note that the statements Lemma~\ref{lem:dualty} also hold if generic identifiability is considered instead. 
	In some circumstances where identifiability, or generic identifiability, of the original network is difficult to analyze, the dual results in Lemma~\ref{lem:dualty} can be applied, leading to a simpler analysis on the basis of the corresponding transpose network. We use the following example to illustrate our point.

	\begin{exam}
		Consider the dynamic network in Fig.~\ref{fig:orginalnet}, where we are investigating identifiability of the SIMO (single input multiple output) subsystem $G_{\mathcal{N}^+_1,1}$. In this subsystem the input $w_1$, and the outputs $w_2$, ${w}_3$ and $w_4$ are excited, while ${w}_2$ and ${w}_5$ are measured.

		The available tool in the current literature to analyze the identifiability of $G_{\mathcal{N}^+_1,1}$ is the disconnecting set-based approach in \cite{shi2020partial}, by which it is required to determine identifiability of each module in $G_{\mathcal{N}^+_1,1}$ separately. For instance,
		to check the identifiability of $G_{41}$ in the original network, according to \cite[Thm.3]{shi2020partial} we need to find an excited vertex set $\bar{\mathcal{R}}$ and a measured vertex set $\bar{\mathcal{C}}$ such that  a disconnecting set $\mathcal{D}$ from $\bar{\mathcal{R}} \cup \mathcal{N}_1^+ \setminus \{4\}$ to $ \bar{\mathcal{C}}$ is found to fulfill the constraints:  $b_{\bar{\mathcal{R}} \rightarrow \mathcal{D}} = |\mathcal{D}|$ and $b_{ \{4\}\cup\mathcal{D} \rightarrow \bar{\mathcal{C}}} = |\mathcal{D}| + 1$. In this example, we find  $\bar{\mathcal{R}} = \{1, 2, 3\}$, $\bar{\mathcal{C}} = \{2,5,6,7\}$, and $\mathcal{D} = \{3,4,5\}$ for (generic) identifiability of $G_{41}$.

		On the other hand, it follows from Lemma~\ref{lem:dualty} that identifiability of $G_{\mathcal{N}^+_1,1}$ in the original network is equivalent to identifiability of $\tilde{G}_{1,\mathcal{N}^-_1}$ in its transpose network shown in Fig.~\ref{fig:transposenet}. $\tilde{G}_{1,\mathcal{N}^-_1}$ is a MISO subsystem with $\tilde{w}_2$, $\tilde{w}_3$, $\tilde{w}_4$, and $\tilde{w}_5$ as the inputs and $\tilde{w}_1$ as the output, while from the network graph, it can be observed that:
		\begin{equation} \label{eq:exT1r}
			[\tilde{T}]_{1,\tilde{\mathcal{R}}} = \tilde{G}_{1, \mathcal{N}_1^-} [\tilde{T}]_{\mathcal{N}_1^-,\tilde{\mathcal{R}}} ,
		\end{equation}
		where $\tilde{\mathcal{R}} = \{2, 5, 6, 7\}$ and $\mathcal{N}_1^- = \{2, 3, 4, 5\}$. The analysis of this MISO subsystem actually requires a simpler procedure than the analysis in \cite{shi2020partial}. 
		
		Although vertex 5 in $\mathcal{N}_1^-$ is not measured, the mapping from $\tilde{\mathcal{R}}$ to vertex $5$ in the transpose network can be obtained as follows. 
		We take the second row of the equation $(I - \tilde{G}) \tilde{T} \tilde{R} = \tilde{R}$, and permute it as 
		\begin{equation*}
			\begin{bmatrix}
				- \tilde{G}_{25} & 1 & 0_{1 \times 3}
			\end{bmatrix}
			\begin{bmatrix}
				[\tilde{T}]_{5,\tilde{\mathcal{R}}} \\ [\tilde{T}]_{2,\tilde{\mathcal{R}}} \\ \star
			\end{bmatrix}
			= [\tilde{R}]_{2,\star} = \mathbf{e}_2,
		\end{equation*}
		with $\mathbf{e}_2$ the second column vector in the identity matrix. The transfer function $[\tilde{T}]_{5,\tilde{\mathcal{R}}}$ is then represented as 
		$[\tilde{T}]_{5,\tilde{\mathcal{R}}} = \tilde{G}_{25}^{-1}(\mathbf{e}_2 - [\tilde{T}]_{2,\tilde{\mathcal{R}}}) = [\tilde{T}]_{2,5}^{-1}(\mathbf{e}_2 - [\tilde{T}]_{2,\tilde{\mathcal{R}}})$.
		Therefore, the transfer matrix $[\tilde{T}]_{\mathcal{N}_1^-,\tilde{\mathcal{R}}}$ in \eqref{eq:exT1r} can be rewritten as 
		\begin{equation*}
			[\tilde{T}]_{\mathcal{N}_1^-,\tilde{\mathcal{R}}} = \begin{bmatrix}
				[\tilde{T}]_{\{3,4\},\tilde{\mathcal{R}}} \\
				[\tilde{T}]_{2,5}^{-1} (\mathbf{e}_2 - [\tilde{T}]_{2,\tilde{\mathcal{R}}})
			\end{bmatrix},
		\end{equation*}
		where all the elements can be acquired from measurement data. By applying the result in \cite{henk2018identifiability}, we know $[\tilde{T}]_{\mathcal{N}_1^-,\tilde{\mathcal{R}}}$ is invertible,  since there exists a so-called \textit{constrained set} of four vertex disjoint paths from $\tilde{\mathcal{R}}$ to $\mathcal{N}_1^-$. As a result, $\tilde{G}_{1,\mathcal{N}^-_1}$ can be uniquely obtained by multiplying both sides of \eqref{eq:exT1r} by $[\tilde{T}]_{\mathcal{N}_1^-,\tilde{\mathcal{R}}}^{-1}$. With Lemma~\ref{lem:dualty} this also implies identifiability of ${G}_{\mathcal{N}^+_1,1}$ in the original network model set.
		
		\label{ex:1}
		\begin{figure}[t] 
			\begin{minipage}[t]{0.5\linewidth}
				\centering
				\includegraphics[width=0.95\textwidth]{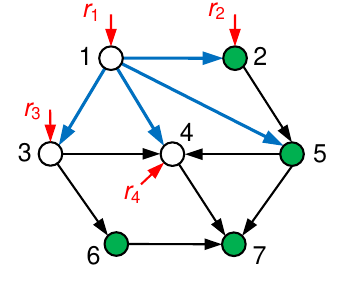}
				\subcaption{Original Network}
				\label{fig:orginalnet}
			\end{minipage}%
			\begin{minipage}[t]{0.5\linewidth}
				\centering
				\includegraphics[width=0.98\textwidth]{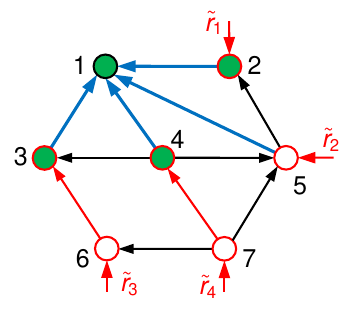}
				\subcaption{Transpose Network}
				\label{fig:transposenet}
			\end{minipage}
			\caption{A dynamic network (a) and its transpose network (b). The filled green vertices are measured, while the excited vertices are indicated by the red arrows. Identifiability of $G_{41}$ in the original network model set is equivalent to that of $\tilde{G}_{14}$ in the transpose network model set.}
			\label{fig:OneModuleCutNew}
		\end{figure}
	\end{exam}
	
	%\subsection{Branchings in Dynamic Networks}
	%
	%\begin{definition}[Branching Truncation]
	%	A branching in a directed graph $\mathcal{G}$ is the maximal subgraph in which every vertex and all its descendants have exact one in-neighbour. A graph $\mathcal{G}_t$ is branching truncation of $\mathcal{G}$ if it is obtained by removing all branchings in $\mathcal{G}$.
	%\end{definition}
	%
	%
	%\begin{pro} 
	%	Let $\mathcal{G}_t$ is branching truncation of $\mathcal{G}$. Define a network model associating with $\mathcal{G}_t$ where the root of each branching is measured if and only there is a vertex in the branching is measured. Then, the branching-truncated network is identifiable if and only if the subnetwork of the original network is identifiable.
	%\end{pro}
	
	\section{Identifiability Analysis in Acyclic Networks}
	\label{sec:Acyclic}
	
	In this section, the identifiability problem is investigated for dynamic networks with acyclic topology. We present two sufficient conditions for identifiability in acyclic networks, where the first condition provides an analysis approach to check identifiability and generic identifiability, and the second condition leads to a synthesis approach for selecting locations of sensors and actuators to achieve network identifiability.
	
	In most of the existing studies e.g., \cite{hendrickx2018identifiability,henk2018identifiability,weerts2018identifiability,shi2020IFAC}, necessary conditions for identifiability of dynamic networks have been provided in the setting of either full excitation or full measurement. A necessary condition for generic identifiability in the case of partial excitation and measurement can be found in \cite{Bazanella2019CDC}, which states that a dynamic network is generically
	identifiable only if every vertex in the network is either excited or measured. Furthermore, the following necessary condition has been provided in \cite{cheng2021necessary}.
	
	\begin{pro}
		\label{pro:necessary1}
		Consider the network model set $\mathcal{M}$ in \eqref{eq:modelset} with $\mathcal{R} \subseteq \mathcal{V}$ and $\mathcal{C} \subseteq \mathcal{V}$ the excited and measured vertices. If $\mathcal{M}$ is identifiable from $[T]_{\mathcal{C},\mathcal{R}}$, then
		\begin{subequations} \label{eq:rankcond}
			\begin{align}
				\rank \left([T]_{\mathcal{N}_j^-,\mathcal{R}}(q, \theta)\right) = |\mathcal{N}_j^-|, 
				\label{eq:rankcond1}
				\\ \text{and} \ 
				\rank \left([T]_{\mathcal{C},\mathcal{N}_j^+}(q, \theta)\right) = |\mathcal{N}_j^+|,
				\label{eq:rankcond2}
			\end{align}
		\end{subequations}
		hold for each $j \in \mathcal{V}$ and for all $\theta \in \Theta$.
		%	, where $\mathcal{P}_j^-$ and $\mathcal{P}_j^+$ are defined as
		%	\begin{align}
			%		\mathcal{P}_j^-: = \{i \in \mathcal{N}_j^- \mid G_{ji}(\theta)~\text{is parameterized in}~\mathcal{M} \}, \label{eq:Pj-}\\
			%		\mathcal{P}_j^+: = \{i \in \mathcal{N}_j^+ \mid G_{ij}(\theta)~\text{is parameterized in}~\mathcal{M} \}. \label{eq:Pj+}
			%	\end{align}
	\end{pro}
	%\begin{proof}
	%	The necessity of \eqref{eq:rankcond1} can be proved following a similar reasoning as in Theorem~2 in  \cite{weerts2018identifiability} for the full measurement case. Then, the necessity of \eqref{eq:rankcond2} has been shown in \cite{hendrickx2018identifiability,henk2018identifiability}. 
	%	The necessity of \eqref{eq:rankcond2}
	%	can also be validated using
	%	 Lemma~\ref{lem:dualty}. Let $\mathcal{G}$ be the graph consistent with the dynamic network \eqref{eq:net} and $\mathcal{G}^\prime$ be the transpose graph of $\mathcal{G}$. Note that identifiability of $\mathcal{M}$ is equivalent to identifiability of $\mathcal{M}^\prime$, which is the model set of its transpose network. Therefore, we have $	\rank ([\tilde{T}]_{\tilde{\mathcal{N}}_i^-,\tilde{\mathcal{R}}}) = |\tilde{\mathcal{N}}_i^-|$ in $\mathcal{M}^\prime$, where $\tilde{\mathcal{R}}$ and $\tilde{\mathcal{N}}_i^-$ are exited vertices and in-neighbors of vertex $i$ in $\mathcal{G}^\prime$. As a result, \eqref{eq:rankcond2} holds, which is obtained from the duality relations: $\tilde{\mathcal{R}} = \mathcal{C}$, $\tilde{\mathcal{N}}_i^- = \mathcal{N}_i^+$ and $\tilde{T} = T^\top$ implied in Definition~\ref{defn:transpose}. 
	%		$\hfill\blacksquare$ 
	%\end{proof}
	
	Note that the generic rank of a transfer matrix can be characterized using a graphical concept, called vertex-disjoint paths \cite{hendrickx2018identifiability,weerts2018single}. Specifically, for any submatrix $T_{\mathcal{B},\mathcal{A}}$ of the parameterized transfer matrix $T(q,\theta)$ in \eqref{eq:T},
	\begin{equation} \label{eq:rankvdp}
		\mathrm{gRank}\left([T]_{\mathcal{B},\mathcal{A}}(q, \theta)\right) = b_{\mathcal{A} \rightarrow \mathcal{B}},
	\end{equation}
	where $b_{\mathcal{A} \rightarrow \mathcal{B}}$ is the maximum number of vertex-disjoint paths from  $\mathcal{A} \subseteq \mathcal{V}$ to $\mathcal{B} \subseteq \mathcal{V}$.
	With the relation \eqref{eq:rankvdp}, the two rank conditions in \eqref{eq:rankcond} can be reinterpreted based on the underlying graph of the model set: If a network model set $\mathcal{M}$ in \eqref{eq:modelset} is (generically)  identifiable, then 
	\begin{align}
		b_{\mathcal{R} \rightarrow \mathcal{N}_j^-}  = |\mathcal{N}_j^-|, 
		\ \text{and} \ 
		b_{\mathcal{N}_j^+ \rightarrow \mathcal{C}}  = |\mathcal{N}_j^+|,
		\label{eq:rankcond2graph}
	\end{align}
	for each vertex $j \in \mathcal{V}$.
	
	In the case that the vertices are either all excited or all measured, the conditions in \eqref{eq:rankcond} and \eqref{eq:rankcond2graph} become sufficient for identifiability and generic identifiability, respectively, see the results in e.g., \cite{hendrickx2018identifiability,weerts2018identifiability,henk2018identifiability}. However, these conditions are not sufficient to determine identifiability of networks where not all vertices are excited or all vertices are measured.

	\subsection{Sufficient Condition for Identifiability in Acyclic Networks}
	
	By exploiting the structural property of acyclic graphs, we present sufficient conditions for the identifiability analysis in dynamic networks with acyclic topology.
	
	Notice that acyclic graphs do not contain directed cycles, and thus any path between two distinct vertices has a finite length.  {\blue A path $p$ from $i_{n}$ to $i_1$ in an acyclic graph is a sequence of edges $(i_n, i_{n-1})$,  $(i_{n-1}, i_{n-2})$, ...,  $(i_{3}, i_{2}),$ $(i_{2}, i_1)$, where each node $i_k$, $k = 1, 2, ...,n$, appears once on the path. We then denote the transfer function of the path $p$ as $T^p_{i_n \rightarrow i_1}: = G_{i_1i_2} G_{i_2i_3} \cdots G_{i_{n-2}i_{n-1}}G_{i_{n-1}i_n}$.} The following property of the transfer matrix $T$ of an acyclic network is then given. 
	
	%Moreover, we adopt the reachability concept in graph theory. In a directed graph, a vertex $j$ is \textit{reachable} from $i$ if there is a path from $i$ to $j$. Based on that, we define the following reachability sets.
	%\begin{definition}
	%	In a directed acyclic graph, the set of vertices that can be reached from a vertex set $\mathcal{X}$ is called the \textbf{reachable set} of $\mathcal{X}$, denoted by $\mathscr{R}_\mathcal{X}$. Furthermore, a \textbf{constrained reachable set} of $\mathcal{X}$ with regard to $\mathcal{Y}$, denoted by  $\mathscr{R}_{\mathcal{Y},\mathcal{X}}$, is a set of vertices that are reachable from $\mathcal{X}$ via $\mathcal{Y}$. 
	%\end{definition}
	%
	%It is remarked that if  $\mathscr{R}_{\mathcal{Y},\mathcal{X}}$ is a  constrained reachable set of $\mathcal{X}$ with regard to $\mathcal{Y}$, then $\mathcal{Y}$ is a disconnecting set from $\mathcal{X}$ to $\mathscr{R}_{\mathcal{Y},\mathcal{X}}$.

	\begin{lemma} \label{lem:Tji}
		Consider the network model in \eqref{eq:net} with the underlying acyclic graph $\mathcal{G}$. Then, the following statements hold.
		\begin{enumerate}
			\item $T_{ii} = 1$, for all vertices in $\mathcal{G}$;
			
			\item $T_{ji} = 0$, if $j \notin \mathscr{R}_i$;
			
			\item {\blue $T_{ji} = \sum_{p \in \mathcal{P}} T^p_{i \rightarrow j}$, for all $j \in \mathscr{R}_i$},
		\end{enumerate}
		where the set $\mathscr{R}_i$ collects all the descendants vertex $i$, and {\blue $\mathcal{P}$ is the set of all paths from $i$ to $j$}. 
	\end{lemma}

	For given sets of excited and measured vertices, identifiability of $G$ essentially reflects if we can uniquely obtain $G$ from the submatrix  $[T]_{\mathcal{C},\mathcal{R}}$. Lemma~\ref{lem:Tji} indicates that every nonzero entry $T_{ji}$ in the matrix $[T]_{\mathcal{C},\mathcal{R}}$ corresponds to a set of directed paths from $i$ to $j$. Based on that, we present an iterative procedure to check the identifiability of individual modules in an acyclic network, where we say a  path is \textit{unknown} at the $k$-th iteration if it contains at least one parameterized edge (module) whose identifiability has not been determined in the previous $k-1$ iterations.
	
	Initially, all the paths containing parameterized edges are unknown, and we present Proposition~\ref{pro:TjiGvu} to iteratively determine the identifiability of each parameterized module in the network. Modules that are shown to be identifiable after the current iteration will be considered as \textit{known} in subsequent iterations.

	\begin{pro} \label{pro:TjiGvu}
		Consider the network model set $\mathcal{M}$ in \eqref{eq:modelset} with $\mathcal{R} \cup \mathcal{C} = \mathcal{V}$ and all the models satisfying Assumption~\ref{Assum}. 
		Then, $G_{\mu \nu}$ is identifiable in $\mathcal{M}$  from $T_{\mathcal{C},\mathcal{R}}$ at the $k$-th iteration, if there exists an excited vertex $i$ and a measured vertex $j$ such that every unknown path from $i$ to $j$ at the $k$-th iteration contains the directed edge $(\nu, \mu)$.
	\end{pro}
	% 	\textcolor{red}{Alternative text for Proposition:} \\
	% \textcolor{red}{Consider the network model set $\mathcal{M}$ in \eqref{eq:modelset} with $\mathcal{R} \cup \mathcal{C} = \mathcal{V}$ and all the models satisfying Assumption~\ref{Assum}. 
		% 		Then, $G_{\mu \nu}$ is identifiable in $\mathcal{M}$  from $T_{\mathcal{C},\mathcal{R}}$ at the $k$-th iteration, if there exists an excited vertex $i$ and a measured vertex $j$ such that every unknown path from $i$ to $j$ at the $k$-th iteration contains the directed edge $(\nu, \mu)$.}

	% \medskip 
	The detailed proof of Proposition~\ref{pro:TjiGvu} is provided in Appendix~\ref{ap:pro:TjiGvu}.  When implementing the condition in Proposition~\ref{pro:TjiGvu}, we may encounter one of the following two cases. One, we find an excited vertex $i$ and a measured vertex $j$ such that all the paths from $i$ to $j$ pass through a common edge $(\nu,\mu)$, then the parameterized module $G_{\mu  \nu }$ is identifiable and becomes a known module in subsequent iterations. Two, the paths from an excited vertex $i$ to a measured vertex $j$ can be divided into two sets: one set only contains known paths (i.e. all the modules on the paths are known), and the other set $\mathcal{P}$ contains unknown paths and each of them passes through a common edge $(\nu,\mu)$. Then the parameterized module $G_{\mu  \nu }$ is identifiable according to Proposition~\ref{pro:TjiGvu}.

	The following corollary presents some special results of Proposition~\ref{pro:TjiGvu}. 
	\begin{coro} \label{coro:singlepath}
		Consider the network model set $\mathcal{M}$ in \eqref{eq:modelset} with $\mathcal{R} \cup \mathcal{C} = \mathcal{V}$ and all the models satisfying Assumption~\ref{Assum}. Then the following statements hold.
		\begin{enumerate}
			\item If there exists a unique directed path from an excited vertex $i$ to a measured vertex $j$, then all the modules in the path are identifiable in $\mathcal{M}$ from $T_{ji}$. 
			\item A tree network is identifiable if and only if the root is excited, and all the leaves are measured \cite{Bazanella2019CDC}.
		\end{enumerate}
	\end{coro}
	The first statement is obtained directly from Proposition~\ref{pro:TjiGvu}. The second statement has been presented in \cite{Bazanella2019CDC}, while it can be proved alternatively in an easy way by using the first statement, as there is always a unique directed path from the root to every leaf in a tree network. For some acyclic networks, the results in Proposition~\ref{pro:TjiGvu} and Corollary~\ref{coro:singlepath} provide an efficient analysis tool for identifiability of single modules or even a full network.  
	\begin{exam}
		\begin{figure}[!tp]\centering
			\includegraphics[scale=0.9]{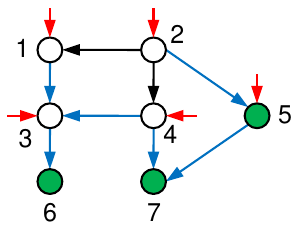}	
			\caption{An acyclic network, in which identifiability of all the modules can be verified by using Proposition~\ref{pro:TjiGvu}.}
			\label{fig:singlepath}
		\end{figure}
		With this example, we demonstrate the use of Proposition~\ref{pro:TjiGvu} and Corollary~\ref{coro:singlepath} in the acyclic network shown in Fig.~\ref{fig:singlepath}, where vertices in $\mathcal{R} = \{1,2,3,4, 5\}$ are excited, and vertices in $\mathcal{C} = \{5,6,7\}$ are measured. We consider the identifiability from the submatrix $[T]_{\mathcal{C},\mathcal{R}}$. 
		
		Observe that $T_{61}$ represents a unique path from vertex $1$ to vertex $6$, then both $G_{31}$ and $G_{63}$ are identifiable in  $\mathcal{M}$ due to condition (a) in Corollary~\ref{coro:singlepath}.
		Similarly, the mappings $T_{52}$, $T_{64}$, $T_{74}$, $T_{75}$ represent unique paths between corresponding vertices, which leads to identifiability of modules $G_{52}$, $G_{34}$, $G_{74}$, and $G_{75}$, respectively. 
		
		Now, all the modules indicated by the blue color are identifiable and thus become known for determining identifiability of $G_{42}$ and $G_{12}$. Since the edge $(2,4)$ appears in the only unknown path from vertex $2$ to vertex $7$, it follows from Proposition~\ref{pro:TjiGvu} that $G_{42}$ is identifiable in $\mathcal{M}$. 
		Then, the module $G_{12}$ is the only unknown module in $\mathcal{M}$, which is also the only unknown edge in the paths from  $2$ to $6$ and hence identifiable. As a result, the entire network is identifiable from $[T]_{\mathcal{C},\mathcal{R}}$.
	\end{exam}
	
	The conditions in Proposition~\ref{pro:TjiGvu} and Corollary~\ref{coro:singlepath} require finding a path between a specific pair of vertices for checking identifiability in each iteration. In some cases, e.g., tree networks or the network in Fig.~\ref{fig:singlepath}, these conditions are able to determine  identifiability of the full network, i.e. all the modules in a network. However, for networks with general topology, the conditions may only be able to determine the identifiability of a subset of modules in the network, rather than the full network. For instance, if we modify the network in Fig.~\ref{fig:singlepath} by adding a new edge $(4,6)$, then we cannot determine the identifiability of modules $G_{3,4}$ and $G_{6,4}$ in this modified network by using the conditions in Proposition~\ref{pro:TjiGvu} and Corollary~\ref{coro:singlepath}.
	
	Next, we develop a new condition to check identifiability of the model set of an acyclic network in the general case. The following main result is provided to determine identifiability based on a vertex-wise inspection.
	
	\begin{theorem}
		\label{thm:identifiability}  
		Consider the network model set $\mathcal{M}$ in \eqref{eq:modelset} with all the models satisfying Assumption~\ref{Assum}. Let $\mathcal{R} \subseteq \mathcal{V}$ and $\mathcal{C} \subseteq \mathcal{V}$ be the sets of excited and measured vertices such that $\mathcal{R} \cup \mathcal{C} = \mathcal{V}$.  The model set $\mathcal{M}$ is identifiable from $[T]_{\mathcal{C},\mathcal{R}}$ if the following conditions hold:
		\begin{enumerate}
			
			\item Every measured vertex $j$ in the network satisfies
			\begin{equation} \label{eq:thm:rank0}
				\rank\left([T]_{\mathcal{N}_j^-,\mathcal{R}}\right) = |\mathcal{N}_j^-|;
			\end{equation}
			
			\item For each unmeasured (and excited) vertex $j$, there exist a set of measured vertices $\mathcal{C}_j$ and a set of excited vertices $\mathcal{R}_j$ such that $\widehat{\mathcal{N}}^-_j \subset  \mathcal{C}_j$, $j \in \mathcal{R}_j$, and
			\begin{align}
				&\rank \left([T]_{\mathcal{C}_j,\mathcal{R}_j}\right) = |\mathcal{R}_j|, 
				\label{eq:thm:rank1}
				\\
				&  \rank\left([T]_{\mathcal{C}_j,(\mathcal{R}_j \cup \mathcal{S}_j) \setminus j}\right) = |\mathcal{R}_j| - 1,
				%		&\rank\left([T]_{\mathcal{C}_j,\mathcal{R}_j \setminus j}\right) = \rank\left([T]_{\mathcal{C}_j,(\mathcal{R}_j \cup \mathcal{S}_j) \setminus j}\right),
				\label{eq:thm:rank2}
			\end{align}
			where $\widehat{\mathcal{N}}^-_j \subseteq {\mathcal{N}}^-_j$ collects all the unexcited in-neighbors of $j$, and 
			$\mathcal{S}_j : = \widehat{\mathcal{N}}^-_j \cup \left(\bigcup_{i\in {\mathcal{N}}^-_j}^{} \mathcal{N}^+_i\right)$. 
		\end{enumerate}
	\end{theorem}
	
	In Theorem~\ref{thm:identifiability}, we leave out the dependence on the delay operator $q$ and parameters $\theta$ in the transfer matrices for ease of notation, while the equations \eqref{eq:thm:rank0}, \eqref{eq:thm:rank1}, and \eqref{eq:thm:rank2} should hold for all $\theta \in \Theta$. The proof is given in Appendix~\ref{ap:thm:identifiability}. 
	
	The conditions in Theorem~\ref{thm:identifiability} are now discussed. 
	The two conditions provide vertex-wise checks for measured and unmeasured vertices in the network separately. Condition (a) requires a sufficient number of excitation signals to the in-neighbors of each measured vertex, which coincides with the rank condition for network identifiability in the full-measurement case \cite{weerts2018identifiability}.  
	Condition (b), on the other hand, is presented for checking the identifiability of the MISO subsystem $\mathcal{G}_{\mathcal{N}_j^+,j}$ regarding each unmeasured vertex $j$ in the network. For an unmeasured vertex $j$, the set $\mathcal{S}_j$ is fixed, which is a union of the $j$'s unexcited in-neighbors (i.e. $\widehat{\mathcal{N}}^-_j$) and the out-neighbors of all the $j$'s in-neighbors,  (i.e. $\bigcup_{i\in {\mathcal{N}}^-_j}^{} \mathcal{N}^+_i$).  Then, to check condition (b), we find a set of excited vertices $\mathcal{R}_j$ containing $j$ and a set of measured vertices $\mathcal{C}_j$ that includes all the unexcited in-neighbors of $j$, and further check if $\mathcal{R}_j$ and $\mathcal{C}_j$ satisfy the two rank requirements in \eqref{eq:thm:rank1} and \eqref{eq:thm:rank2}. 
	
\begin{figure}[!tp]\centering
		\centering
		\includegraphics[scale=0.35]{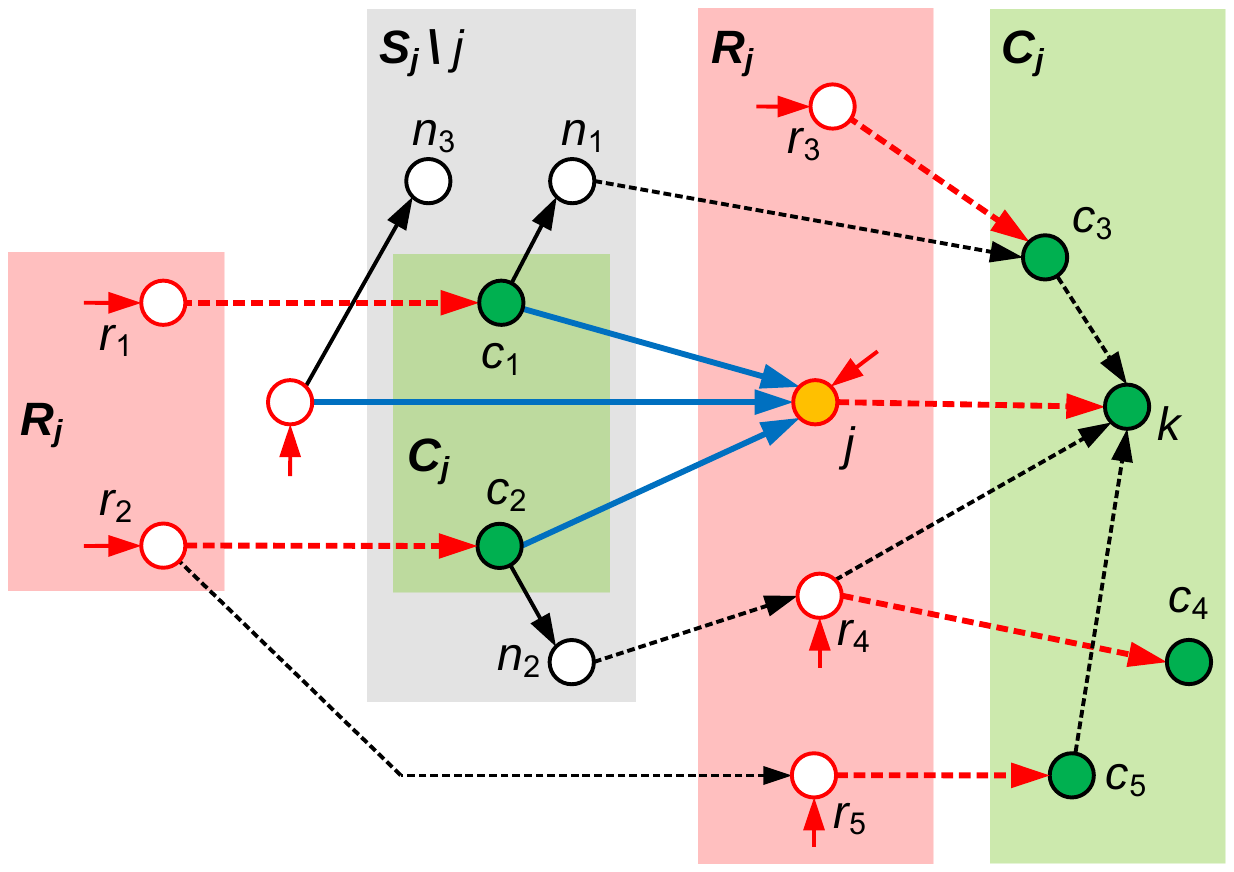}	
		\caption{An illustration of the second condition in Theorem~\ref{thm:identifiability}, where the green circles indicate measured vertices, and dashed arrows represent paths.}
		\label{fig:theorem2a}
	\end{figure}

	Next, we illustrate the two constraints in \eqref{eq:thm:rank1} and \eqref{eq:thm:rank2} by using Fig.~\ref{fig:theorem2a}. 
	% where $\mathcal{R}_j = \{ j, r_1, r_2, r_3, r_4, r_5\}$ and $\mathcal{C}_j = \{k, c_1, c_2, c_3, c_4, c_5\}$. From \cite{bazanella2017identifiability,hendrickx2018identifiability}, the full column rank of $[T]_{\mathcal{C}_j,\mathcal{R}_j}$ in \eqref{eq:thm:rank1} implies that there are $|\mathcal{R}_j| = 6$ vertex-disjoint paths from $\mathcal{R}_j$ to $\mathcal{C}_j$, indicated by the red dashed arrows, in which there is a directed path from $j$ to a measured vertex $k \in \mathcal{C}_j$. Besides, the columns of matrix $[T]_{\mathcal{C}_j,(\mathcal{R}_j \cup \mathcal{S}_j) \setminus j}$ in \eqref{eq:thm:rank2} contains two compartments: $[T]_{\mathcal{C}_j,\mathcal{R}_j \setminus j}$ and $[T]_{\mathcal{C}_j,\mathcal{S}_j \setminus j}$, 
	% with $\mathcal{S}_j = \{j, c_1, c_2, n_1, n_2, n_3 \}$ collecting all the unexcited in-neighbors of $j$, i.e. $\widehat{\mathcal{N}}_j^- = \{c_1, c_2\}$, and all the out-neighbors of $j$'s in-neighbors $\bigcup_{i\in {\mathcal{N}}^-_j}^{} \mathcal{N}^+_i = \{j, n_1, n_2, n_3\}$. From  \eqref{eq:thm:rank1} and \eqref{eq:thm:rank2}, we have
	% \begin{align} \label{eq:thm:rank12}
		% 	\rank \left([T]_{\mathcal{C}_j,\mathcal{R}_j\setminus j}\right) &  =  |\mathcal{R}_j| - 1
		% 	\\
		% 	& = \rank \left(\begin{bmatrix}
			% 		[T]_{\mathcal{C}_j,\mathcal{R}_j \setminus j} & [T]_{\mathcal{C}_j,\mathcal{S}_j \setminus j}
			% 	\end{bmatrix}\right). \nonumber
		% \end{align}
	Here, the set $\mathcal{S}_j = \{j, c_1, c_2, n_1, n_2, n_3 \}$ collects all the unexcited in-neighbors of $j$, i.e. $\widehat{\mathcal{N}}_j^- = \{c_1, c_2\}$, and all the out-neighbors of $j$'s in-neighbors $\bigcup_{i\in {\mathcal{N}}^-_j}^{} \mathcal{N}^+_i = \{j, n_1, n_2, n_3\}$. To apply the two conditions \eqref{eq:thm:rank1} and \eqref{eq:thm:rank2}, we find a set of excited vertices $\mathcal{R}_j$ with $j \in \mathcal{R}_j$ and a set of measured vertices $\mathcal{C}_j$ that includes $\widehat{\mathcal{N}}_j^- \subseteq \mathcal{C}_j$ and a measured descendant $k$. In Fig.~\ref{fig:theorem2a}, we choose  
	$\mathcal{R}_j = \{ j, r_1, r_2, r_3, r_4, r_5\}$ and $\mathcal{C}_j = \{k, c_1, c_2, c_3, c_4, c_5\}$ such that the transfer matrix $[T]_{\mathcal{C}_j,\mathcal{R}_j}$ in \eqref{eq:thm:rank1} has full column rank. This is guaranteed by the $|\mathcal{R}_j| = 6$ vertex-disjoint paths from $\mathcal{R}_j$ to $\mathcal{C}_j$, indicated by the red dashed arrows, in which there is a directed path from $j$ to a measured vertex $k \in \mathcal{C}_j$.
	Besides, with $(\mathcal{R}_j \cup \mathcal{S}_j)\setminus j = \{j, r_1, r_2, r_3, r_4, r_5, c_1, c_2, n_1, n_2, n_3 \}$, $\rank \left([T]_{\mathcal{C}_j,\mathcal{R}_j\setminus j}\right)  =  |\mathcal{R}_j| - 1$, as there are 5 vertex disjoint paths from $(\mathcal{R}_j \cup \mathcal{S}_j)\setminus j$ to $\mathcal{C}_j$. Therefore, \eqref{eq:thm:rank2} is fulfilled. 
	
	\begin{remark}
		\label{rem:disconnectset}
		Note that the columns of matrix $[T]_{\mathcal{C}_j,(\mathcal{R}_j \cup \mathcal{S}_j) \setminus j}$ in \eqref{eq:thm:rank2} contain two compartments: $[T]_{\mathcal{C}_j,\mathcal{R}_j \setminus j}$ and $[T]_{\mathcal{C}_j,\mathcal{S}_j \setminus j}$, and 
		\begin{align} \label{eq:thm:rank12}
			\rank \left([T]_{\mathcal{C}_j,\mathcal{R}_j\setminus j}\right)  =  
			\rank \left(\begin{bmatrix}
				[T]_{\mathcal{C}_j,\mathcal{R}_j \setminus j} & [T]_{\mathcal{C}_j,\mathcal{S}_j \setminus j}
			\end{bmatrix}\right).  
		\end{align}%
		It follows from \eqref{eq:rankvdp} that there are no vertex-disjoint paths from $\mathcal{S}_j \setminus j$ to $\mathcal{C}_j$ in the graph that are vertex-disjoint with the ones from $\mathcal{R}_j \setminus j$ to $\mathcal{C}_j$. Further, \eqref{eq:thm:rank2} can be interpreted using the concept of disconnecting set, as follows. If \eqref{eq:thm:rank12} holds, there is a disconnecting set $\mathcal{D}$ with $|\mathcal{D}| = |\mathcal{R}_j| - 1$, whose removal will lead to absence of paths from $\mathcal{R}_j \setminus j$ to $\mathcal{C}_j$ and from $\mathcal{S}_j \setminus j$ to $\mathcal{C}_j$, simultaneously. In Fig.~\ref{fig:theorem2a}, this set can be chosen as $\mathcal{D} = \{c_1, c_2, c_3,r_4,r_5\}$. Note that the unexcited in-neighbors of $j$, i.e. $\widehat{\mathcal{N}}_j^-$, are contained in $\mathcal{C}_j$, $\mathcal{S}_j$, and $\mathcal{D}$ simultaneously.
	\end{remark}

	The rank conditions in Theorem~\ref{thm:identifiability} can also be represented by means of vertex disjoint paths for generic identifiability of $\mathcal{M}$ as in e.g., \cite{bazanella2017identifiability,hendrickx2018identifiability,weerts2018single,shi2020partial}. Specifically, the following graph-based result is obtained directly from  Theorem~\ref{thm:identifiability}.
	
	\begin{coro} \label{coro:paths}
		Consider the network model set $\mathcal{M}$ in \eqref{eq:modelset} with $\mathcal{R}$ and $\mathcal{C}$ the sets of excited and measured vertices such that $\mathcal{R} \cup \mathcal{C} = \mathcal{V}$. The model set $\mathcal{M}$ is generically identifiable from $[T]_{\mathcal{C},\mathcal{R}}$ if the following conditions hold: 
		\begin{enumerate}
			\item Every measured vertex $j$ in the network satisfies
			\begin{equation} \label{eq:rem:path0}
				b_{\mathcal{R} \rightarrow \mathcal{N}_j^-} = |\mathcal{N}_j^-|;
			\end{equation}
			
			\item For every   unmeasured but excited vertex $j$, there exist a set of measured vertices $\mathcal{C}_j$ and a set of excited vertices $\mathcal{R}_j$ such that $\widehat{\mathcal{N}}^-_j \subset  \mathcal{C}_j$, $j \in \mathcal{R}_j$, and
			\begin{align}
				& b_{\mathcal{R}_j \rightarrow \mathcal{C}_j}  = |\mathcal{R}_j|, 
				\label{eq:rem:path1}\\
				&b_{(\mathcal{R}_j \cup \mathcal{S}_j) \setminus j \rightarrow \mathcal{C}_j} = |\mathcal{R}_j| - 1,
				\label{eq:rem:path2}
			\end{align}
			where $\widehat{\mathcal{N}}^-_j \subseteq {\mathcal{N}}^-_j$ collects all the unexcited in-neighbors of $j$, and 
			$\mathcal{S}_j : = \widehat{\mathcal{N}}^-_j \cup \left(\bigcup_{i\in {\mathcal{N}}^-_j}^{} \mathcal{N}^+_i\right)$. 
			
		\end{enumerate}
	\end{coro}	
	
	{\blue The conditions of Corollary~\ref{coro:paths} or Theorem~\ref{thm:identifiability} lead to the following result, with the proof contained in the proof of Theorem~\ref{thm:identifiability}.
		\begin{pro}
			Consider a network model set $\mathcal{M}$ with the same setup as in Theorem~\ref{thm:identifiability}. If the two conditions in Theorem~\ref{thm:identifiability} or Corollary~\ref{coro:paths} hold, then all the sources in the network are excited, and all the sinks are measured. 
		\end{pro}
	}
	
	This is actually necessary for the identifiability of a network model set in the partial excitation and measurement setting, for more details see \cite{bazanella2017identifiability,Bazanella2019CDC}. Moreover, we can show that if condition (b) in Corollary~\ref{coro:paths} is satisfied, then the following result holds, the proof of which is added in Appendix~\ref{ap:pro:thm-necess}.
	\begin{pro}\label{pro:thm-necess}
		Consider the setting of Corollary~\ref{coro:paths}. If \eqref{eq:rem:path1} and \eqref{eq:rem:path2} hold, then
		\begin{align}
			b_{\mathcal{R}\rightarrow \mathcal{N}_j^-} & = |\mathcal{N}_j^-|, ~\forall~j \in \mathcal{R}, j \notin \mathcal{C}. 
			\label{eq:pro:cond1}
			%		\\
			%		b_{\mathcal{N}_i^+ \rightarrow \mathcal{C}}  &= |\mathcal{N}_i^+|, ~\forall~ i \in \mathcal{C} \cup \mathcal{R}
			%		\label{eq:pro:cond2}
		\end{align}
	\end{pro}
	%If there is vertex $i \in \mathcal{V}$ such that \eqref{eq:pro:cond2} is not satisfied. 
	The condition \eqref{eq:pro:cond1} holds for all unmeasured and excited vertices in the network, and combining this with \eqref{eq:rem:path0} leads to the {\blue first necessary condition} in \eqref{eq:rankcond2graph}. {\blue However, inferring the second necessary condition from Corollary~\ref{coro:paths} is not straightforward, since the proof of the theorem is based on analyzing the in-neighbors of each vertex rather than the out-neighbors. }

	In the following example, we illustrate how to determine identifiability on the basis of network topology by using Corollary~\ref{coro:paths}.
	\begin{exam}
		This example demonstrates how to use the conditions in Corollary~\ref{coro:paths} to check network identifiability of the model set $\mathcal{M}$ of the full acyclic network in Fig.~\ref{fig:acyclicnet}. First, it is not hard to verify that each measured vertex in $\mathcal{C} = \{3, 4, 5\}$ satisfies \eqref{eq:rem:path0}, with $\mathcal{R} = \{1,2,4,6\}$.  Then, to verify generic identifiability of $\mathcal{M}$, we apply condition (b) in Corollary~\ref{coro:paths} to check the unmeasured vertex $1$, $2$, and $6$.
		
		Note that for the source vertices $j=1,6$, we have $\mathcal{N}_j^{-} = \emptyset$, which implies $\mathcal{S}_j = \emptyset$. Thus, the condition (b) in Corollary~\ref{coro:paths} becomes: $b_{\mathcal{R}_j \rightarrow \mathcal{C}_j}  = |\mathcal{R}_j|$ and $b_{\mathcal{R}_j  \setminus j \rightarrow \mathcal{C}_j} = |\mathcal{R}_j| - 1$, which is satisfied by simply choosing $\mathcal{R}_j = \{j\}$ and $\mathcal{C}_j = \{5\}$. Next, vertex $2$ is checked, which has the in-neighbors $\mathcal{N}_2^- = \{1 ,5\}$ and $\widehat{\mathcal{N}}_2^- = \{5\}$ unexcited. Therefore, we have 
		$\mathcal{S}_2 =   \widehat{\mathcal{N}}^-_2 \cup \left(\bigcup_{i\in {\mathcal{N}}^-_2}^{} \mathcal{N}^+_i\right) = \{2, 4, 5\}$. 
		Consider $\mathcal{R}_2 = \{2,4\}$ and $\mathcal{C}_2 = \{3,5\}$ that satisfy \eqref{eq:rem:path1}, i.e. $b_{\mathcal{R}_j \rightarrow \mathcal{C}_j} = 2$. It can be verified that
		$ b_{(\mathcal{R}_2 \cup \mathcal{S}_2) \setminus 2 \rightarrow \mathcal{C}_2} = b_{\{4,5\} \rightarrow \{3,5\}} = 1$, i.e. condition (b) in Corollary~\ref{coro:paths} is also fulfilled.
		As a result, all the vertices in this network satisfy the conditions in Corollary~\ref{coro:paths}, and thus the model set $\mathcal{M}$ is generically identifiable from $[T]_{\mathcal{C},\mathcal{R}}$. 
		
		\begin{figure}[!tp]\centering
			\includegraphics[scale=0.9]{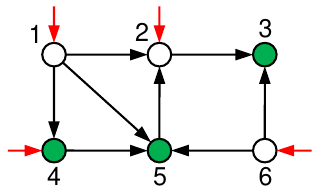}	
			\caption{An acyclic network, in which all the vertices satisfy the conditions in Theorem~\ref{thm:identifiability}. Thus, the network model set is identifiable.}
			\label{fig:acyclicnet}
		\end{figure}
		
	\end{exam}

	To prove the result in Theorem~\ref{thm:identifiability}, we utilize the hierarchical structure in an acyclic network that partitions the vertices into different layers. We refer to Appendix~\ref{ap:thm:identifiability} for the detailed reasoning. The conditions in Theorem~\ref{thm:identifiability} are the result of a vertex-wise analysis starting from the top layer to the lower ones. It is worth noting that \eqref{eq:thm:rank0} is also a necessary condition for measured vertices, while for every unmeasured vertex, condition (b) in Theorem~\ref{thm:identifiability} needs to be satisfied to guarantee that each module with this vertex as an output is identifiable in $\mathcal{M}$. In contrast to  Theorem~\ref{thm:identifiability}, we can also apply dual reasoning that starts from the measured bottom layer in an acyclic network. Then, we will find an alternative condition for identifiability as follows.  
	\begin{coro}
		\label{coro:identifiability_dual}  
		Consider the network model set $\mathcal{M}$ in \eqref{eq:modelset}. Let $\mathcal{R} \subseteq \mathcal{V}$ and $\mathcal{C} \subseteq \mathcal{V}$ be the sets of excited and measured vertices such that $\mathcal{R} \cup \mathcal{C} = \mathcal{V}$.  The model set $\mathcal{M}$ is identifiable from $[T]_{\mathcal{C},\mathcal{R}}$ if the following conditions hold:
		\begin{enumerate}
			\item Every excited vertex $i$ in the network satisfies
			\begin{equation*}
				\rank\left([T]_{\mathcal{C},\mathcal{N}_i^+}\right) = |\mathcal{N}_i^+|;
			\end{equation*}
			
			\item For each unexcited but measured vertex $i$, there exist sets of measured vertices $\mathcal{C}_i$ and excited vertices $\mathcal{R}_i$ such that  $\widehat{\mathcal{N}}^+_i \subset  \mathcal{R}_i$, $i \in \mathcal{C}_j$, and
			\begin{align*}
				&\rank \left([T]_{\mathcal{C}_i,\mathcal{R}_i}\right) = |\mathcal{C}_i|, \\
				& \rank\left([T]_{(\mathcal{C}_i \cup \mathcal{S}_i) \setminus i, \mathcal{R}_i}\right) = |\mathcal{C}_i| -1,
			\end{align*}
			where $\widehat{\mathcal{N}}^+_i \subseteq {\mathcal{N}}^+_i$ collects all the unmeasured out-neighbors of $i$, and 
			$\mathcal{S}_i : = \widehat{\mathcal{N}}^+_i \cup \left(\bigcup_{j\in {\mathcal{N}}^+_i}^{} \mathcal{N}^-_j\right)$. 
		\end{enumerate}
	\end{coro}
	The results can be proved by directly applying Lemma~\ref{lem:dualty}, considering the transpose network of $\mathcal{M}$, and thus is omitted here. Similarly, we can extend the rank conditions in Corollary~\ref{coro:identifiability_dual} to path-based conditions for generic identifiability of $\mathcal{M}$, where, dual to condition (b) in Corollary~\ref{coro:paths}, we will have $b_{\mathcal{N}_i^+ \rightarrow \mathcal{C}} = |\mathcal{C}_i|$ and $b_{\mathcal{R}_i \rightarrow \mathcal{C}_i} = b_{\mathcal{R}_i \rightarrow (\mathcal{C}_i \cup \mathcal{S}_i) \setminus i} + 1 =  |\mathcal{C}_i|$ replacing  the rank conditions in Corollary~\ref{coro:identifiability_dual}(b).

	\begin{remark}
		All the rank conditions in Theorem~\ref{thm:identifiability} and Corollary~\ref{coro:identifiability_dual} can be regarded as a generalization of the conditions in the full measurement case \cite{weerts2018identifiability} or the full excitation case \cite{hendrickx2018identifiability,henk2018identifiability}.
		If all the vertices are measured,  only \eqref{eq:thm:rank0} is required. When all the vertices are excited, we turn to verify condition (b) of Theorem~\ref{thm:identifiability} or Corollary~\ref{coro:identifiability_dual} for each vertex $i$. Moreover, we can reformulate \eqref{eq:thm:rank1} and \eqref{eq:thm:rank2} by considering the identifiability of the modules $G_{\mathcal{N}_i^+,i}$ in the MISO subsystem. To this end, $j$ in \eqref{eq:thm:rank2} is replaced by $\mathcal{N}_i^+$ so that  \eqref{eq:thm:rank2} holds immediately. Now, we have $\rank([T]_{\mathcal{C}_j,\mathcal{R}_j}) = |\mathcal{R}_j|$ for some sets $\mathcal{C}_j \subseteq \mathcal{C}$ and $\mathcal{R}_j \subseteq \mathcal{R}$, where $\mathcal{N}_i^+ \subseteq \mathcal{R}_j$. Clearly, it is equivalent to $\rank([T]_{\mathcal{C},\mathcal{N}_i^+}) = |\mathcal{N}_i^+|$, that is the necessary and sufficient condition for network identifiability in the full excitation case.
	\end{remark}
	
	%Although these condition is interpreted completely using the analysis based on graphs, and it is instrumental for developing a graph-theoretical approach to allocate excitations and measurements in acyclic networks.

	\subsection{Allocation of Actuators and Sensors in Acyclic Networks}
	
	In this part we present an alternative result for verifying identifiability of an acyclic network. While the rank conditions in Theorem~\ref{thm:identifiability} and Corollary~\ref{coro:identifiability_dual} have provided an instrumental tool for verifying identifiability of an acyclic network model set, the concerned conditions are formulated in terms of a vertex-wise analysis. That is to say, for each vertex, it has to be checked separately whether the conditions for identifiability are satisfied. Besides the computational burden that such a procedure entails, it is also an analysis result that is less suitable for addressing the synthesis question: which nodes to excite and which nodes to measure so as to guarantee network identifiability?  This question is about generating a valid Excitation and Measurement Pattern (EMP), see 
	\cite{mapurunga2021optimal} for the definition. In the subsequent part, we present a synthesis scheme to allocate excitation and measurement signals for identifiability of an acyclic network model set. This will lead to alternative and more compact identifiability conditions. 
	
	%
	%	In this part, we present a synthesis scheme to allocate excitation and measurement signals for identifiability of an acyclic network model set.
	%The rank condition in Theorem~\ref{thm:identifiability} and Corollary~\ref{coro:identifiability_dual} have provided an instrumental tool for verifying identifiability of an acyclic network model set. However, these conditions are formulated in terms of  a vertex-wise analysis, which is less effective for addressing the  synthesis problem: which nodes to excite and which nodes to measure so as to guarantee network identifiability? Thus, it is desired to have a more compact identifiability condition to address this problem.
	
	To this end, we resort to a graph covering approach as introduced in \cite{cheng2019allocation} for general (cyclic) dynamic networks. In \cite{cheng2019allocation},  all the vertices in a directed network are measured, and the allocation of excitations is based on decomposing the underlying graph of the network into edge-disjoint \textit{pseudotrees}, whose roots are supposed to be excited for generic identifiability. In this work, we extend the work in  \cite{cheng2019allocation} to the partial excitation and measurement case. Because we are handling acyclic networks here, the disjoint pseudotrees can be relaxed to disjoint directed trees, as explained next.

	A \textit{directed tree}, denoted by $\mathcal{T}$, is a special acyclic graph containing a unique \textit{root} vertex, from which there is exactly one directed path to every other vertex in $\mathcal{T}$. The sinks in $\mathcal{T}$ are also called \textit{leaves}, and the vertices that are neither the root nor leaves of $\mathcal{T}$ are \textit{internal vertices}. Analogously, an anti-tree $\breve{\mathcal{T}}$ is defined, where each vertex in $\breve{\mathcal{T}}$ has exactly one directed path to a unique \textit{root} vertex, while all the sources in $\breve{\mathcal{T}}$ are called \textit{leaves}. 
	
	Then, we define the concept of disjoint trees.
	\begin{definition}[Disjoint trees] \label{def:disjtree}
		In a directed acyclic graph $\mathcal{G}$, two trees $\mathcal{T}_1$ and $\mathcal{T}_2$ are called \textbf{disjoint} if
		\begin{enumerate}
			\item $\mathcal{T}_1$ and $\mathcal{T}_2$ do not share edges, and
			\item {\blue all the edges incident from a vertex are included in the same tree.}
		\end{enumerate}  
	\end{definition}
	
	In parallel with disjoint trees, we can define disjoint anti-trees. Two anti-trees $\mathcal{T}_1$ and $\mathcal{T}_2$ are \textit{disjoint} if $\mathcal{T}_1$ and $\mathcal{T}_2$ do not share edges, and the edges incident to a common vertex are included in the same anti-tree. 
	
	It should be emphasized that there always exists a set of disjoint trees/anti-trees covering all the edges of an acyclic graph, and the covering is not unique. This above statement can be proved directly following \cite{cheng2019allocation}. Based on that, we present a new identifiability condition for acyclic networks as follows.
	
	\begin{theorem} [Tree/Anti-Tree Covering]\label{thm:treecover}
		Consider a network model set $\mathcal{M}$ associating with an acyclic graph $\mathcal{G}$ with $\mathcal{V} = \mathcal{R} \cup \mathcal{C}$. Then, $\mathcal{M}$ is generically identifiable if either of the following two conditions holds:
		\begin{enumerate}
			\item $\mathcal{G}$ can be decomposed into a set of disjoint trees, and for each tree, its root is excited and all the leaves are measured.
			
			\item  $\mathcal{G}$ can be decomposed into a set of disjoint anti-trees, and for each anti-tree, its root is measured, and all the leaves are excited.
		\end{enumerate} 
	\end{theorem}
	
	{\blue Compared to Corollary~1 in \cite{cheng2019allocation}, where all the vertices are assumed to be measured and the goal is to allocate actuators only, Theorem~\ref{thm:treecover} presents generic identifiability conditions for dynamic networks with partial excitations and partial measurements, which can be used to allocate both actuators and sensors.}
		The detailed proof is presented in Appendix~\ref{ap:thm:treecover}. Note that the second condition in Theorem~\ref{thm:treecover} can be proven using the results in Section~\ref{sec:TransposeNet} on transpose networks, but
	it is worth noting that the two conditions in Theorem~\ref{thm:treecover} are dual but not equivalent. Even if there does not exist a disjoint tree covering obeying condition (a), we may still find a set of disjoint anti-trees that covers all the edges of $\mathcal{G}$ and satisfies condition (a). Thus, the two conditions are actually complementary to each other.
	The following example is given to demonstrate this point.

	\begin{exam}
		\begin{figure}
			\centering
			\begin{minipage}[t]{0.5\linewidth}
				\centering
				\includegraphics[width=\textwidth]{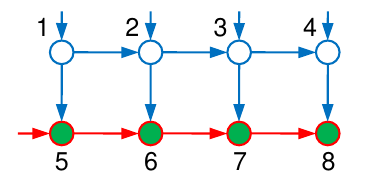}
				\subcaption{}
				\label{fig:treecovering}
			\end{minipage}%
			\begin{minipage}[t]{0.5\linewidth}
				\centering
				\includegraphics[width=\textwidth]{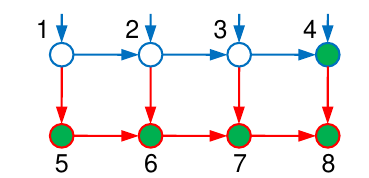}
				\subcaption{}
				\label{fig:antitreecovering}
			\end{minipage}
			\caption{Generic identifiability analysis of an acyclic network using Theorem~\ref{thm:treecover}. (a) The network is decomposed into two disjoint trees, indicated by different colors, and each tree has an excited root and measured leaves. (b) An anti-tree covering is found in a network with the same topology as in (a), where each anti-tree has a measured root and excited leaves.}
		\end{figure}
		Consider a model set $\mathcal{M}$ governed by the dynamic network in Fig.~\ref{fig:treecovering}. It can be found that there are two disjoint trees highlighted by the blue and red colors, covering all the edges of the graph, while their roots (vertices $1$ and $5$) are excited, and the leaf (vertex $8$) is measured. Therefore, it follows from condition~(1) in Theorem~\ref{thm:treecover} that $\mathcal{M}$ is generically identifiable. Note that each of the vertices $2$, $3$, and $4$ can be either excited or measured here, since {\blue they} are neither the roots nor the roots of the blue tree. However, we cannot apply condition (b) to check the generic identifiability of the dynamic network in Fig.~\ref{fig:treecovering}, as we cannot find a set of disjoint anti-trees with measured roots and excited leaves in this graph.
		In contrast, condition (b) can be used to check the generic identifiability of the network in Fig.~\ref{fig:antitreecovering}, where vertex $4$ is measured and vertex $5$ is not excited. Then, we can find two disjoint anti-trees indicated by the blue and red colors, and both anti-trees have measured roots and excited leaves, so that condition~(b) can be applied to yield generic identifiability of $\mathcal{M}$ in this case. 
	\end{exam}

	Compared to the sufficient conditions in Theorem~\ref{thm:identifiability}, the conditions in Theorem~\ref{thm:treecover} are more conservative.
	% in the sense of identifiability analysis. 
	This conservativeness is embodied in the following facts:
	\begin{enumerate}
		\item If condition (a) of Theorem~\ref{thm:treecover} holds, then $|\mathcal{N}_j^-| \leq 1$ for every unmeasured vertex $j$ in the network;
		\item If condition (b) of Theorem~\ref{thm:treecover} holds, then $|\mathcal{N}_i^+| \leq 1$ for every unexcited vertex $i$ in the network.
	\end{enumerate}
	The above necessary conditions can be verified by using the disjointness property of trees and anti-trees. They are 
	implicitly required when applying the two conditions in Theorem~\ref{thm:treecover}.
	Despite the conservativeness compared to Theorem~\ref{thm:identifiability}, Theorem~\ref{thm:treecover} provides more compact conditions on the level of a full graph rather than a vertex-wise analysis. This will facilitate the development of a synthesis procedure for the signal allocation problem in acyclic networks.
	
	\subsection{Algorithm}

	Given an acyclic network without any excitations and measurements, {\blue the objective of the synthesis problem is} to \textit{allocate a minimum number of actuators and sensors, i.e., $|\mathcal{R}| + |\mathcal{C}|$, to achieve generic identifiability} of the network.
	
\begin{figure}[t] 
		\begin{minipage}[t]{.5\textwidth}
			\centering 
			\includegraphics[width=0.82\textwidth]{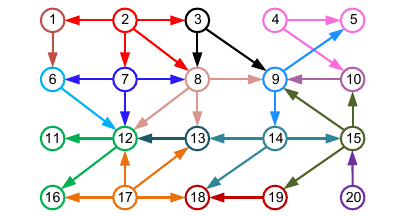}
			\subcaption{Initialization}
			\label{fig:covering1}
		\end{minipage} 
		\begin{minipage}[t]{.5\textwidth} 
			\centering
			\includegraphics[width=0.8\textwidth]{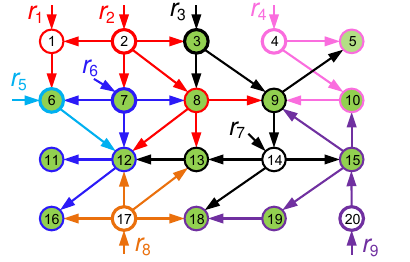}
			\subcaption{Tree Covering and Signal Allocation}
			\label{fig:covering2}
		\end{minipage}
		\caption{Illustration of the tree covering method for actuator/sensor allocation. (a) An acyclic graph is covered by a set of disjoint minimal trees, which are indicated by different colors. (b) A tree covering with a minimum cardinality. The roots and leaves of the trees are to be excited and measured, respectively. Furthermore, {\blue $r_5$ is redundant} for generic identifiability by applying Corollary~\ref{coro:paths} and thus can be removed.}
		\label{fig:covering}
	\end{figure}

	{\blue Generally, minimizing $|\mathcal{R}|+|\mathcal{C}|$ is a  difficult problem, particularly in the context of large-scale networks where an optimal solution is often challenging to find. Therefore, we devise a graph-theoretical algorithm on the basis of Theorem~\ref{thm:treecover} to tackle this problem in a heuristic manner.}
	The main steps follow similarly as in \cite{cheng2019allocation}, where we have presented a heuristic algorithm to find a set of disjoint pseudotrees to cover all the edges in a cyclic graph. In the following, we devise a two-step scheme for signal allocation in acyclic networks.
	The detailed procedure is elaborated, with an illustration of an example in Fig.~\ref{fig:covering}.

	(1) \textit{Tree Covering.} 
	For any acyclic graph, we can always find a set of disjoint trees covering all its edges. Due to the identifiability condition in Theorem~\ref{thm:treecover}, it is desirable to have a covering with a minimum number of trees. However, this minimal covering problem is a   combinatorial optimization problem, whose solution is not unique. Therefore, we adopt a heuristic algorithm to find a locally optimal solution to the tree covering problem, see Algorithm~\ref{alg:covering}.
	
	\begin{algorithm}[t]
		\caption{Tree Covering of Acyclic Graphs}
		\hspace*{\algorithmicindent} \textbf{Input}: An acyclic graph $\mathcal{G}$ 
		\begin{algorithmic}[1]
			\STATE Decompose $\mathcal{G}$ into a set of minimal trees.
			\REPEAT
			\STATE 
			Merge two trees into a single tree if their union remains a tree. 
			\UNTIL There are no mergeable trees.
		\end{algorithmic}
		\label{alg:covering}
	\end{algorithm}

	In Algorithm~\ref{alg:covering}, a minimal tree is a tree that consists of a vertex and all its out-neighbors with all the edges incident from the vertex to the out-neighbors. Note that any acyclic graph admits such a decomposition and any two minimal trees that {\blue have no} common vertices are disjoint \cite{cheng2019allocation}. In Fig.~\ref{fig:covering1}, a graph of 20 vertices is partitioned into 16 minimal trees indicated by different colors. Starting from the minimal tree covering, the merging step will reduce the total number of trees required to cover the graph. 
	%A tree $\mathcal{T}_1$ can be merged into $\mathcal{T}_2$  whose root is the root of $\mathcal{T}_2$. 
	%In Fig.~\ref{fig:covering1}, the tree rooted at 1 can be merged into the tree rooted at 2. However, the trees rooted at 2 and 3 cannot be merged, since their combination is no longer a tree as vertices 2 and 3 share the same out-neighbor 8. 
	For the graph in Fig.~\ref{fig:covering1}, a tree covering is generated, shown in Fig.~\ref{fig:covering2}, which contains only 7 disjoint trees. 
	It can be verified that for this specific example, the resulting tree covering has a minimum cardinality, {\blue i.e. the number of trees in the covering is minimal.}
	
	\begin{remark}
		Note that the result of the tree covering for a given graph is dependent on the order of merging. The minimum covering is not always guaranteed, while the algorithm always converges to feasible tree covering. 
		% a local minimum can be achieved. 
		In \cite{cheng2019allocation}, a particular order of merging is proposed to yield an effective solution to the minimum covering problem. The details, therefore, are omitted in this paper. 
	\end{remark}
	%In the first stage, we merge a tree $\mathcal{T}_1$ into $\mathcal{T}_2$, if $\mathcal{T}_1$ can only be merged into $\mathcal{T}_2$. For instance, the tree rooted at 19 in Fig.~\ref{fig:covering1} can be merged into the unique tree rooted at 15, thus we can merge them first. While the tree rooted at 9 can be merged into either the tree rooted at 3 or 8, we then leave them to be dealt with later. In the second stage, we randomly merge the remaining pairs until there are no more trees that can be merged. 

	(2) \textit{Signal Allocation.}
	A direct application of condition (a) in Theorem~\ref{thm:treecover} implies that we can allocate excitation signals at the roots of the trees in the resulting covering and measure all the leaves of these trees, while the rest of the vertices in the graph can be either excited or measured. condition (a) in Theorem~\ref{thm:treecover} then guarantees the generic identifiability of the model set of this acyclic network.

	Note that if two or more trees have a vertex in common, this shared vertex has to be a leaf of a tree and hence is required to be measured, which may lead to a set of measured roots in the trees, coinciding with $\mathcal{R} \cap \mathcal{C}$, while the other vertices in $\mathcal{G}$ are either excited or measured. For the network in Fig.~\ref{fig:covering2}, the roots $\mathscr{R}_t = \{2, 3, 4, 6, 7, 17, 20\}$ are to be excited, and the leaves $\mathscr{L}_t  = \{3, 5, 6, 7, 8, 9, 10, 11, 12, 13, 15, 16, 18\}$ are to be measured. The remaining vertices $1$, $14$, and $19$ are internal vertices in the trees that can be either excited or measured. In this example, we choose to excite vertices $1$ and $14$ and measure vertex $19$. Then, the network model set is generically identifiable due to Theorem~\ref{thm:treecover}.

	Observe that when $|\mathcal{R}| + |\mathcal{C}|$ is minimized, it equivalently yields a minimum number of vertices that are excited and measured simultaneously, i.e. $|\mathcal{R} \cap \mathcal{C}|$ is minimized. {\blue It is worth pointing out that finding a covering with a minimal number of trees does not necessarily imply that $|\mathcal{R}|+|\mathcal{C}|$ is minimized.}
	To further reduce the number of required signals, a subsequent operation is taken to check the redundant excitation signals assigned to the vertices in the set $\mathcal{R} \cap \mathcal{C}$ by using Corollary~\ref{coro:paths}. Removing these unnecessary excitation signals further gives a smaller cardinality of $\mathcal{R} \cap \mathcal{C}$.
	
	Specifically, we inspect each tree $\mathcal{T}_k$ whose root is an element in $\mathcal{R} \cap \mathcal{C}$ that is excited and measured simultaneously. If the network model set remains generically identifiable after removing the excitation signal on the root of $\mathcal{T}_k$, i.e. all measured vertices in $\mathcal{T}_k$ still satisfy \eqref{eq:rem:path0}, and all the excited vertices in $\mathcal{T}_k$ fulfill the conditions in \eqref{eq:rem:path1} and \eqref{eq:rem:path2}, then
	we can take out the excitation signal allocated at the root of $\mathcal{T}_k$.

	We consider the network in Fig.~\ref{fig:covering2} as an example. Observe that $\mathcal{R} \cap \mathcal{C} = \{3, 6, 7\}$, which are all measured root vertices. If we remove $r_5$ on vertex $6$, the resulting network will remain 
	generically identifiable since $b_{\mathcal{R}\setminus \{6\} \rightarrow \mathcal{N}_k^{-}} = |\mathcal{N}_k^{-}|$ holds for both measured vertex $k=7$ and $8$ in the cyan tree. Note that we cannot remove $r_3$ in the black tree, since vertex $9$ in this tree has four in-neighbors and according to condition~(a) of Corollary~\ref{coro:paths}, four vertex-disjoint paths are needed from the excited vertices to its in-neighbors, which would not be satisfied if $r_3$ is removed  {\blue Following a similar reasoning, we cannot eliminate $r_6$ by analyzing the number of in-neighbors of vertex
	$8$ in the blue tree.}

	%\begin{algorithm}[t]
	%	\caption{Allocation of Actuators and Sensors in Acyclic Networks}
	%	\begin{algorithmic}[1]
		%		\STATE  Fina a minimum number of disjoint trees $\mathcal{T}_1, \mathcal{T}_2, \cdots, \mathcal{T}_\kappa$ that cover all the edges of the acyclic graph $\mathcal{G}$.
		%		
		%		\STATE For each tree $\mathcal{T}_i$, place the root, denoted by $\tau (\mathcal{T}_i)$ in $\mathcal{R}$, and all the leaves of $\mathcal{T}_i$ in $\mathcal{C}$. 
		%		
		%		\IF{a vertex in $\mathcal{G}$ is shared by two or more trees} 
		%		\STATE Put this vertex in $\mathcal{C}$;
		%		\ELSE
		%		\STATE Put this vertex in either $\mathcal{C}$ or $\mathcal{R}$.
		%		\ENDIF
		%		
		%		\FOR{$i = 1, ..., \kappa$, and $\tau (\mathcal{T}_i) \in \mathcal{C}$, $\tau (\mathcal{T}_i) \in \mathcal{R}$} 
		%		\STATE Removing $\tau (\mathcal{T}_i)$ from $\mathcal{R}$.
		%		\IF{The conditions in Theorem~\ref{thm:identifiability} does not hold} 
		%		\STATE Add $\tau (\mathcal{T}_i)$ back into $\mathcal{R}$.
		%		\ENDIF
		%		\ENDFOR
		%		
		%		\RETURN $\mathcal{R}$ and $\mathcal{C}$ as actuator and sensor vertices
		%	\end{algorithmic}
	%	\label{alg:treecovering}
	%\end{algorithm}

	%\begin{exam}
	%	We use a ten-vertex acyclic network in Fig.~\ref{fig:covering} to illustrate the procedure of Algorithm~\ref{alg:treecovering}. First, a disjoint tree covering of the graph is found with a minimum cardinality, and these trees are distinguished by different colors in Fig.~\ref{fig:covering}. Then, we select their roots and leaves to be excited and measured, respectively. Thus, Next, we check if the number of excitation signals can be further reduced.
	%	

	%\end{exam}
	
\begin{figure}[t] 
		\centering
		\includegraphics[width=0.4\textwidth]{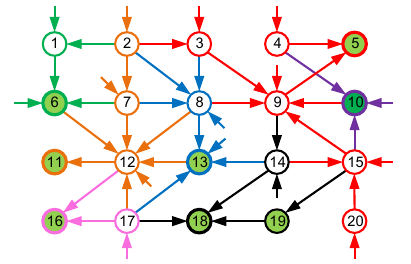}
		\caption{Illustration of the anti-tree covering method for actuator/sensor allocation. The acyclic graph in Fig.~\ref{fig:covering} is now covered by a set of disjoint anti-trees, which are indicated by different colors. In each anti-tree, the root is measured, and the leaves are excited. Furthermore, it can be verified that the measurement at {\blue vertex 10} is redundant for generic identifiability and thus can be removed.}
		\label{fig:covering-dual}
	\end{figure}

	\begin{remark}
		The above procedure utilizes the first condition in Theorem~\ref{thm:treecover} to find the two potential sets $\mathcal{R}$ and $\mathcal{C}$ for allocating actuators and sensors, respectively. Alternatively, we can also design a different scheme according to the second condition in Theorem~\ref{thm:treecover}, where an anti-tree covering of $\mathcal{G}$ is considered. Particularly, we first seek for a minimum set of disjoint anti-trees that cover all the edges of $\mathcal{G}$ and then allocate sensors and actuators at the roots and leaves of the anti-trees, respectively. To pursue a minimum number of measured vertices, we apply Corollary~\ref{coro:identifiability_dual} to check if the vertices in $\mathcal{R} \cap \mathcal{C}$ are necessary to be measured. We apply the anti-tree covering scheme to the network in Fig.~\ref{fig:covering} and obtain a set of anti-trees as shown in Fig.~\ref{fig:covering-dual}, where the roots of the anti-trees, i.e. $\mathscr{R}_t = \{5, 6, 10, 11, 13, 16, 18, 19\}$, are measured, and the all the leaves in $\mathscr{L}_t = \{2,4,7,8,10,12,13,14,15,17,20\}$ are excited. The internal vertices, including 1, 9, and 19, can be either excited or measured. Thereby, we have {\blue $\mathcal{R} \cap \mathcal{C} = \{6, 10, 13\}$. Notice that $b_{\mathcal{N}_k^{+} \rightarrow \mathcal{C}\setminus \{10\}  } = |\mathcal{N}_k^{+}|$ holds for all the vertices in the anti-tree rooted at vertex 10. Therefore, we can further remove vertex 10 from $\mathcal{C}$, and the network model set will remain 
		generically identifiable.}
	\end{remark}
	
	\section{Conclusions}
	\label{sec:conclusion}
	
	In this paper, we have analyzed identifiability in acyclic dynamic networks where only partial excitation and measurement signals are available. The concept of transpose networks is introduced, whose identifiability is shown to be equivalent to identifiability of its corresponding original network. We have also presented a necessary condition for identifiability of general networks, where the identifiability problem is interpreted as solving a system of nonlinear equations with the parameterized modules as indeterminate variables.
	
	The other major contribution of this paper is to present two sufficient conditions for identifiability of acyclic networks. The first condition can be implemented to analyze identifiability of a given acyclic network based on a vertex-wise check. The second condition, on the other hand, is established using the concept of disjoint tree/anti-tree covering and is instrumental in developing an algorithmic procedure that selects excitation and measured signals such that an acyclic network is generically identifiable.

	\section*{Appendix}
	\renewcommand{\thesubsection}{\Alph{subsection}}

	\subsection{Proof of Lemma \ref{lem:dualty}}
	\label{ap:lem:dualty}
	
	From Definition~\ref{defn:netid}, a module ${G}_{ij}$ is identifiable in the model set $\mathcal{M}$  if and only if the following implication holds:
	\begin{equation} \label{eq:implication1}
		C T(q,\theta_0) R = CT(q,\theta_1)  R  \Rightarrow {G}_{ij}(q,\theta_0) = {G}_{ij}(q,\theta_1),
	\end{equation}
	for all $\theta_0,\theta_1 \in \Theta$. Analogously, the network identifiability of a module $\tilde{G}_{ji}$ in the model set $\mathcal{M}(\mathcal{G}^\prime)$ is equivalent to the implication
	\begin{multline} \label{eq:implication2}
		R^\top  (I - \tilde{G}(q, \theta_0))^{-1} C^\top = R^\top (I - \tilde{G}(q, \theta_1))^{-1} C^\top  \\ \Rightarrow \tilde{G}_{ji}(q,\theta_0) = \tilde{G}_{ji}(q,\theta_1),
	\end{multline}
	Observe that $ (C T R)^\top = R^\top (I - \tilde{G})^{-1} C^\top$, which is the transfer matrix of the transpose network. Therefore, the equations on the left-hand side of implications \eqref{eq:implication1} and \eqref{eq:implication2} are equivalent, from which the first claim immediately follows. The second statement on the overall model set is obtained by applying the above reasoning to all the modules in the network \eqref{eq:net}.

	\subsection{Proof of Lemma \ref{lem:Tji}}
	\label{ap:lem:Tji}
	
	The analytic function $T$ can be expanded in the Taylor series as 
	\begin{align}
		T = (I - G)^{-1} = I + \sum_{k=1}^{\infty} G^k.
	\end{align}
	To prove the two statements in the lemma, it is sufficient to show that, for any $k \geq 1$, 
	$[G^k]_{ii} = 0$, and $[G^k]_{ij}$ is the sum of all directed paths from $j$ to $i$ of length $k$.
	
	We proceed with the proof by induction on $k$. For $k = 1$, the result holds immediately, since $G^1 = G$, where the main diagonal entries are all $0$, while each nonzero off-diagonal entry $G_{ij}$ presents a path of length $1$ from $j$ to $i$. Assume now that the inductive hypothesis holds for some $k \geq 1$. 
	Note that any path of length $k+1$ from $j$ to $i$, with $i \neq j$, consists of an edge $G_{iz}$ for some vertex $z \in \mathcal{N}_i^-$ and a path of length $k$ from $j$ and $z$. It then gives
	\begin{align} \label{eq:k+1}
		[G^{k+1}]_{ij} = \sum_{z\in \mathcal{N}_i^-}^{ } G_{iz} [G^k]_{zj},
	\end{align}
	where $[G^k]_{zj}$ represents the sum of all directed paths of length $k$ from $j$ to $z$. Therefore, the expression of $[G^{k+1}]_{ij} $ in \eqref{eq:k+1} implies the sum of all directed paths of length $k+1$ from $j$ to $i$.  For the case that $i = j$ in \eqref{eq:k+1}, $[G^k]_{zi} = 0$ holds since $\mathcal{G}$ is acyclic, e.g., there is no directed circle from vertex $i$ to itself, which leads to  $[G^{k+1}]_{ii} = 0$. That completes the proof.

	\subsection{Proof of Proposition \ref{pro:TjiGvu}}
	\label{ap:pro:TjiGvu}
	Let $\mathcal{P}$ be the set of all paths from $i$ to $j$. If the condition in Proposition \ref{pro:TjiGvu} holds, then the mapping $T_{ji}$, without loss of generality,  can be represented as
	\begin{equation}\label{eq:TjiGvu}
		T_{ji} = T_{j \mu} G_{\mu\nu } T_{\nu i} + T_{\mathrm{rem}},
	\end{equation}
	{\blue where $T_{j \mu} = \sum_{p \in \mathcal{P}_1} T^p_{\mu \rightarrow j}$ and $T_{\nu i} = \sum_{p \in \mathcal{P}_2} T^p_{\nu \rightarrow i}$, where the sets $\mathcal{P}_1$ and $\mathcal{P}_2$ collect all the paths from $\mu$ to $j$ and from $i$ to $\nu$, respectively. }$T_{\mathrm{rem}}$ represents all the paths from $i$ to $j$ that exclude the ones in $\mathcal{P}$ and do not contain the edge $(\nu, \mu)$. It implies that all the elements in $T_{\mathrm{rem}}$ are known from the previous $k-1$ iterations. Note that the set $\mathcal{P}$ collects all the unknown paths corresponding to the transfer matrix $T_{j \mu} G_{\mu\nu } T_{\nu i}$, which is known and computed as $ T_{ji} - T_{\mathrm{rem}}$.  
	
	Then, the identifiability of $G_{\mu \nu}$ in $\mathcal{M}$ is analyzed for four situations by considering different excitation/measurement patterns of nodes $\mu$ and $\nu$: either $\mu$ is
	excited and $\nu$ is measured, or the contrary or both are excited or
	both are measured. 
	First, if 
	vertex $\nu$ is measured, and $\mu$ is excited, then $T_{j \mu}$ and $T_{\nu i}$ can be obtained so that $G_{\mu\nu } = T_{j \mu}^{-1}(T_{ji} - T_{\mathrm{rem}}) T_{\nu i}^{-1}$.
	Second, if vertex $\nu$ is excited, and $\mu$ is measured, then
	$G_{\mu\nu }$ can be directly obtained from the mapping  $T_{\mu \nu}$, yielding that all the  paths from $\nu$ to $\mu$ other than $G_{\mu\nu }$ are known from $T_{\mathrm{rem}}$. Third, if both vertices $\mu$ and $\nu$ are excited, the two terms $T_{j \mu}$ and $T_{j \nu} = T_{j \mu} G_{ \mu \nu}$ in \eqref{eq:TjiGvu} can be obtained due to excited vertex $i$ and measured vertex $j$. From the two terms, $G_{\mu \nu}$ can be obtained. For the fourth situation where $\mu$ and $\nu$ are both measured, we simply compute $G_{ji}$ from the two terms $T_{\mu i} = G_{\mu \nu} T_{\nu i}$ and $T_{\mu i}$ in \eqref{eq:TjiGvu}, which shows the identifiability of $G_{\nu \mu}$ in $\mathcal{M}$.

	\subsection{Proof of Theorem \ref{thm:identifiability}}
	\label{ap:thm:identifiability}
	In our previous works \cite{shi2020subnetworks,shi2020partial}, the concept of disconnecting sets has been used as a key enabler for analyzing network identifiability and addressing the excitation and measurement allocation problem within a single module setting.
	In this paper, we also capitalize on the properties of disconnecting sets to prove the identifiability conditions for a full network in Theorem~\ref{thm:identifiability}.  
	
	Before proving  Theorem~\ref{thm:identifiability}, two important properties of disconnecting sets are recapped from \cite{schrijver2003book,shi2020subnetworks,shi2020partial}. First, the concept of disconnecting sets is closely related to vertex disjoint paths, as presented in the well-known Menger’s theorem.
	\begin{theorem}\cite{schrijver2003book}
		\label{thm:Menger}
		Let $\mathcal{X}$ and $\mathcal{Y}$ be two vertex subsets in a directed graph $\mathcal{G}$. Let $\mathcal{D}$ be a minimum disconnecting set from $\mathcal{X}$ to $\mathcal{Y}$. Then,
		\begin{equation}
			|\mathcal{D}| = b_{\mathcal{X} \rightarrow \mathcal{Y}},
		\end{equation}
		where $b_{\mathcal{X} \rightarrow \mathcal{Y}}$ is 
		the maximum number of vertex
		disjoint paths from $\mathcal{X}$ to $\mathcal{Y}$.
	\end{theorem}
	
	It is further shown in \cite{shi2020subnetworks} that we can obtain a factorization of the mapping from the vertex signals from $\mathcal{X}$ to $\mathcal{Y}$ in a dynamic network by means of a disconnecting set $\mathcal{D}$ from $\mathcal{X}$ to $\mathcal{Y}$ in the associated graph. This factorization is instrumental for the proof of Theorem~\ref{thm:identifiability}.
	
	\begin{lemma}\cite{shi2020subnetworks,shi2020partial}
		\label{lem:factorization}
		Consider a network model in \eqref{eq:net} associated with a
		directed graph $\mathcal{G}$, and the transfer matrix $T$ in \eqref{eq:T}. Let $\mathcal{D}$ be a disconnecting set from $\mathcal{X}$ to $\mathcal{Y}$ in $\mathcal{G}$. 
		Then there exist proper transfer matrices $K$ and $Q$ such
		that
		\begin{equation} \label{eq:factorization}
			[T]_{\mathcal{Y},\mathcal{X}} = K [T]_{\mathcal{D},\mathcal{X}} =  [T]_{\mathcal{Y},\mathcal{D}} Q,
		\end{equation}
		where $K$ is left invertible if $\rank([T]_{\mathcal{Y},\mathcal{D}}) = |\mathcal{D}|$, and $Q$ is right invertible if $\rank([T]_{\mathcal{D},\mathcal{X}}) = |\mathcal{D}|$.
	\end{lemma}
	
	Thanks to the properties of disconnecting sets in Theorem~\ref{thm:Menger} and Lemma~\ref{lem:factorization},  we are ready to prove Theorem~\ref{thm:identifiability}.

	The proof relies on a particular property of acyclic graphs, that is a hierarchical structure in terms of vertex reachability. Specifically, for any acyclic graph $\mathcal{G}$, we can partition its vertex set $\mathcal{V}$ into different layers (subsets), denoted by $\mathcal{L}_1, \mathcal{L}_2, \cdots, \mathcal{L}_N$ with $N \leq L$, such that the following conditions hold:
	\begin{enumerate}
		\item The vertices in the same layer cannot reach each other; 
		\item For any two layers $\mathcal{L}_\mu$, $ \mathcal{L}_\nu$ with $\mu < \nu$, there is a vertex in $\mathcal{L}_\mu$ that can reach a vertex in $\mathcal{L}_\nu$, but not vise versa.
	\end{enumerate}
	Note that there may be multiple options to partition the vertices in an acyclic graph $\mathcal{G}$ to form a hierarchical structure. In this proof, we consider a specific one where the top layer $\mathcal{L}_1$ is the set of source vertices, while the bottom layer $\mathcal{L}_N$ collects all the sink vertices. Moreover, we can prove that $\mathcal{L}_1$ is excited and $\mathcal{L}_N$ is measured as follows. 
	
	By contradiction, assume a vertex $j \in \mathcal{L}_N$ is unmeasured, then $j$ is excited due to
	$\mathcal{R} \cup \mathcal{C} = \mathcal{V}$. From condition (b) in Theorem~\ref{thm:identifiability}, it is required that there exists a set of measured vertices $\mathcal{C}_j$ such that \eqref{eq:thm:rank1} holds. Since $j$ has no out-neighbors and $j \in \mathcal{R}_j$, the $j$-th column of the matrix $[T]_{\mathcal{C}_j,\mathcal{R}_j}$ is zero, which contradicts \eqref{eq:thm:rank1} that requires full column rank of $[T]_{\mathcal{C}_j,\mathcal{R}_j}$. Therefore, all the vertices in $ \mathcal{L}_N$ are measured. 
	{\blue Next, we can show by contradiction that all the vertices $k \in \mathcal{L}_1$ are excited. Suppose that $k$ is measured but not excited, and let vertex $j$ to be any out-neighbor of $k$. If $j$ is measured. Then we cannot find a set of excited vertices to satisfy condition (a) in Theorem~\ref{thm:identifiability}, since $k \in {\mathcal{N}}_j^-$ is a source. 		
		If $j$ is unmeasured but excited, then \eqref{eq:thm:rank1} and \eqref{eq:thm:rank2} should be fulfilled, which implies
		\begin{align}
			 	\label{eq:thmproof:rank}
			\rank\left([T]_{\mathcal{C}_j,\mathcal{R}_j \setminus j}\right)  = 
			|\mathcal{R}_j| - 1 = 
			\rank\left([T]_{\mathcal{C}_j,(\mathcal{R}_j \cup \mathcal{S}_j) \setminus j}\right).
		\end{align}
		 Note that $k \in \widehat{\mathcal{N}}_j^-$, which means that $k$ is an element in both $\mathcal{S}_j \setminus j$ and $\mathcal{C}_j$. Since $k \notin \mathcal{R}_j\setminus j$ and $k$ is a  source, $[T]_{k,\mathcal{R}_j \setminus j} = 0$, while $[T]_{k,\mathcal{S}_j \setminus j} \ne 0$. Therefore,
		\begin{align*}  
			\rank \left([T]_{\mathcal{C}_j,\mathcal{R}_j\setminus j}\right)  < & 
			\rank \left(\begin{bmatrix}
				[T]_{\mathcal{C}_j,\mathcal{R}_j \setminus j} & [T]_{\mathcal{C}_j,\mathcal{S}_j \setminus j}
			\end{bmatrix}\right)
		 \\ 
		  & =	\rank\left([T]_{\mathcal{C}_j,(\mathcal{R}_j \cup \mathcal{S}_j) \setminus j}\right),  
		\end{align*}%
		which contradicts \eqref{eq:thmproof:rank}. Therefore, all the vertices in $\mathcal{L}_1$ are excited.
     }

% If  hold, then for each vertex $j$ in $\mathcal{G}$, 
%	\begin{equation} \label{eq:thmproofpath}
%		b_{\mathcal{R}\rightarrow \mathcal{N}_j^-} = |\mathcal{N}_j^-|.
%	\end{equation}
%	We refer to \eqref{eq:rem:path0} and   Proposition~\ref{pro:necessary1} for this statement, which are implied by Theorem~\ref{thm:identifiability}. Therefore,  any $k \in \mathcal{L}_1$ should be excited. Otherwise, there will be a out-neighbor $j$  of $k$ that does not satisfy \eqref{eq:thmproofpath}.
%	
%	We can apply an analogous reasoning to the transpose network with graph $\mathcal{G}^\prime$ of the original network, which gives that all the sinks in $\mathcal{G}^\prime$ are measured if the model set of the transpose network is identifiable. With Lemma~\ref{lem:dualty}, it is necessary to have all the sources of $\mathcal{G}$, i.e. the vertices in $\mathcal{L}_1$, are excited for identifiability of $\mathcal{M}$.
	
	Then, identifiability of the network model set $\mathcal{M}$ is equivalent to all the modules $G_{ji}$ being identifiable in  $\mathcal{M}$, for all $i \in \mathcal{N}_j^-$ and $j$ in all the layers. We, therefore, proceed by induction.
	
	First, we analyze the vertices in $\mathcal{L}_2$. For any measured vertex $j \in \mathcal{L}_2$, all its in-neighbors $\mathcal{N}_j^-$ are sources in $\mathcal{L}_1$, which are excited, and $\widehat{\mathcal{N}}^-_j = \emptyset$. From each $k \in \mathcal{N}_j^- \subseteq \mathcal{L}_1$, the directed edge $(i,j)$ is the unique path from $i$ to $j$, thus we obtain from Corollary~\ref{coro:singlepath} that $G_{ji}$ is identifiable in $\mathcal{M}$ for all $i\in \mathcal{N}_j^-$. For the other vertices in $\mathcal{L}_2$ that are not measured, they are supposed to be excited due to $\mathcal{C} \cup \mathcal{R} = \mathcal{V}$. Consider the modules $G_{j, \mathcal{N}_j^-}$ with $j\in \mathcal{L}_2$ and $\mathcal{N}_j^- \subseteq \mathcal{L}_1$ excited. Then, we show identifiability of $G_{j, \mathcal{N}_j^-}$ with both $\mathcal{N}_j^-$ and $j$ excited whenever the condition (b) in Theorem~\ref{thm:identifiability} is satisfied. In the following, we extend our result in \cite{shi2020partial} for single modules to treat the MISO subsystem regarding the unmeasured vertex $j$.
	
%	It is implied from \eqref{eq:thm:rank1} and \eqref{eq:thm:rank2} that 
%	\begin{align*}
%		%	\label{eq:thm:rank3}
%		\rank\left([T]_{\mathcal{C}_j,\mathcal{R}_j \setminus j}\right)  = 
%		|\mathcal{R}_j| - 1 = 
%		\rank\left([T]_{\mathcal{C}_j,(\mathcal{R}_j \cup \mathcal{S}_j) \setminus j}\right).
%	\end{align*}
	Let $\mathcal{D}$ be a minimum disconnecting set from $(\mathcal{R}_j \cup \mathcal{S}_j) \setminus j$ to $\mathcal{C}_j$, where in this case, all the in-neighbors of $j$ are in the layer $\mathcal{L}_1$, which are excited, and hence $\widehat{\mathcal{N}}_j^- = \varnothing$. Then, $\mathcal{S}_j = \bigcup_{i\in {\mathcal{N}}^-_j} \mathcal{N}^+_i$ only contains all the out-neighbors of each in-neighbor of $j$. With the relation between transfer matrix rank and vertex disjoint paths in \cite{vanderWoude1991graph,hendrickx2018identifiability}, the rank equality in \eqref{eq:thmproof:rank} leads to 
	\begin{equation} \label{eq:bbD}
		b_{\mathcal{R}_j \setminus j \rightarrow \mathcal{C}_j}  = b_{(\mathcal{R}_j \cup \mathcal{S}_j) \setminus j \rightarrow \mathcal{C}_j} = |\mathcal{D}|.
	\end{equation}
	It implies that $\mathcal{D}$ is also a disconnecting set from $\mathcal{R}_j \setminus j$ to $\mathcal{C}_j$ as well as a disconnecting set from $\mathcal{N}_i^+ \setminus j$ to $\mathcal{C}_j$ for each $i \in \mathcal{N}^-_j$. Then, we follow  Lemma~\ref{lem:factorization} to obtain that there exist two transfer matrices $Q_1$ and $Q_{2i}$ satisfying
	\begin{subequations}\label{eq:Q2}
		\begin{align} 
			[T]_{\mathcal{C}_j,\mathcal{R}_j \setminus j} &=  [T]_{\mathcal{C}_j,\mathcal{D}} Q_1, 
			\\
			[T]_{\mathcal{C}_j, \mathcal{N}_i^+ \setminus j} 
			&=  [T]_{\mathcal{C}_j, \mathcal{D} } Q_{2i}.
		\end{align}
	\end{subequations}
	Observe that $[T]_{\mathcal{C}_j,\mathcal{R}_j}$ is full column rank due to \eqref{eq:thm:rank1}, which implies that $[T]_{\mathcal{C}_j,\mathcal{R}_j \setminus j}$ is also full column rank, and it is obtained from \eqref{eq:bbD} that $\rank ([T]_{\mathcal{C}_j,\mathcal{R}_j \setminus j}) = |\mathcal{D}|$. Therefore, $Q_1$ is right invertible from Lemma~\ref{lem:factorization} and its right inverse is denoted by $Q_1^\dagger$. Then, it is obtained from \eqref{eq:Q2} that
	\begin{align}\label{eq:subst}
		[T]_{\mathcal{C}_j, \mathcal{N}_i^+ \setminus j} 
		=  [T]_{\mathcal{C}_j,\mathcal{R}_j \setminus j} Q_1^\dagger Q_{2i}.
	\end{align}
	Furthermore, we take the $i$-th column of the equation $C T (I - G) = C$ that can be permuted to 
	\begin{align} 
		\label{eq:ithcolumn} 
		\begin{bmatrix}
			[T]_{\mathcal{C}_j, j} & [T]_{\mathcal{C}_j, \mathcal{N}_i^+ \setminus j} & [T]_{\mathcal{C}_j,i} & \star
		\end{bmatrix}
		\begin{bmatrix}
			-G_{ji} \\   -G_{\mathcal{N}_i^+ \setminus j, i} \\ 1 \\0
		\end{bmatrix}
		=
		C_{\star i},
	\end{align}
	where $C_{\star i}$ denotes the $i$-th column of the matrix $C$, and $C_{\star i}$ is nonzero if vertex $i$ is measured, and $C_{\star i} = 0$ otherwise. 
	Substituting \eqref{eq:subst} into \eqref{eq:ithcolumn} then yields
	\begin{align}\label{eq:TCji}
		%	[T]_{\mathcal{C}_j, i} - C_{\star i} & = [T]_{\mathcal{C}_j, j} G_{ji} + [T]_{\mathcal{C}_j, \mathcal{N}_i^+ \setminus j} G_{\mathcal{N}_i^+ \setminus j, i}
		%	\nonumber\\
		%	& = [T]_{\mathcal{C}_j, j} G_{ji} + [T]_{\mathcal{C}_j,\mathcal{R}_j \setminus j} Q_1^\dagger Q_{2i} G_{\mathcal{N}_i^+ \setminus j, i}.
		[T]_{\mathcal{C}_j, i} - C_{\star i}   = [T]_{\mathcal{C}_j, j} G_{ji} + [T]_{\mathcal{C}_j,\mathcal{R}_j \setminus j} Q_1^\dagger Q_{2i} G_{\mathcal{N}_i^+ \setminus j, i}.
	\end{align}
	Let $[C]_{\star \mathcal{N}_j^-}$ be the submatrix of $C$ that consists of the columns in $C$ indexed by $\mathcal{N}_j^-$. By stacking all the equations in \eqref{eq:TCji} for all $i \in \mathcal{N}_j^-$, we then obtain
	\begin{align}\label{eq:TCjN}
		&[T]_{\mathcal{C}_j, \mathcal{N}_j^-} - [C]_{\star \mathcal{N}_j^-} 
		\\
		&= [T]_{\mathcal{C}_j, j} G_{j,\mathcal{N}_j^-} + [T]_{\mathcal{C}_j,\mathcal{R}_j \setminus j} Q_1^\dagger
		\sum_{i \in \mathcal{N}_j^-}  Q_{2i} G_{\mathcal{N}_i^+ \setminus j, \mathcal{N}_j^-}
		\nonumber\\
		&= 
		\underbrace{\begin{bmatrix}
				[T]_{\mathcal{C}_j, j} & [T]_{\mathcal{C}_j,\mathcal{R}_j \setminus j}
		\end{bmatrix}}_{[T]_{\mathcal{C}_j, \mathcal{R}_j}}
		\begin{bmatrix}
			G_{j,\mathcal{N}_j^-} \\ Q_1^\dagger
			\sum_{i \in \mathcal{N}_j^-}  Q_{2i} G_{\mathcal{N}_i^+ \setminus j, \mathcal{N}_j^-}
		\end{bmatrix},
		\nonumber
	\end{align}
	%	\begin{align}  
		%	[T]_{\mathcal{C}_j, i} - C_{\star i} & = [T]_{\mathcal{C}_j, j} G_{ji} + [T]_{\mathcal{C}_j, \mathcal{N}_i^+ \setminus j} G_{\mathcal{N}_i^+ \setminus j, i}
		%	\nonumber\\ 
		%	& = \underbrace{\begin{bmatrix}
				%		[T]_{\mathcal{C}_j, j} & [T]_{\mathcal{C}_j,\mathcal{R}_j \setminus j}
				%		\end{bmatrix}}_{[T]_{\mathcal{C}_j, \mathcal{R}_j}}
		%	\begin{bmatrix}
			%	G_{ji} \\ Q_1^\dagger Q_2 G_{\mathcal{N}_i^+ \setminus j, i}
			%	\end{bmatrix},
		%	\end{align}
	from which $G_{j,\mathcal{N}_j^-}$ can be uniquely determined from the known terms $[T]_{\mathcal{C}_j, \mathcal{R}_j}$ and $[T]_{\mathcal{C}_j, \mathcal{N}_j^-} - [C]_{\star \mathcal{N}_j^-}$, since  $[T]_{\mathcal{C}_j, \mathcal{R}_j}$ is full column rank as in \eqref{eq:thm:rank1}, and the vertices in $\mathcal{N}_j^- \subseteq \mathcal{L}_1$ are excited. So far, it has been verified that for any vertex $j \in \mathcal{L}_2$ that is measured or excited,  all the modules in $G_{j, \mathcal{N}_j^-}$ are identifiable in the network model set $\mathcal{M}$ if the two conditions in Theorem~\ref{thm:identifiability} hold.

	To proceed with induction on the layers, we assume that the transfer vector $G_{\ell,\mathcal{N}_\ell^-}$ is identifiable in $\mathcal{M}$ for each $\ell \in \widehat{\mathcal{L}} : = \mathcal{L}_1 \cup \mathcal{L}_1 \cup \cdots \cup \mathcal{L}_{k-1}$. We aim to show that $G_{j,\mathcal{N}_j^-}$ is also identifiable, for all  $j \in \mathcal{L}_k$. Notice that the hierarchical structure of the acyclic graph allows for a permutation matrix $P$ such that
	\begin{equation*}
		G = P \begin{bmatrix}
			G_{\widehat{\mathcal{L}},\widehat{\mathcal{L}}} & 0 \\
			G_{\mathcal{V}\setminus\widehat{\mathcal{L}},\widehat{\mathcal{L}}} & G_{\mathcal{V}\setminus\widehat{\mathcal{L}},\mathcal{V}\setminus\widehat{\mathcal{L}}} 
		\end{bmatrix} P^{\top},
	\end{equation*}
	in which all the nonzero entries in $G_{\widehat{\mathcal{L}},\widehat{\mathcal{L}}} $ are identifiable in $\mathcal{M}$. Based on that, we have 
	\begin{equation*}
		T = P^\top (I - P G P^\top)^{-1} P  = P^\top \begin{bmatrix}
			(I - G_{\widehat{\mathcal{L}},\widehat{\mathcal{L}}})^{-1} & 0 \\
			\star & \star
		\end{bmatrix} P,
	\end{equation*}
	and thus the transfer matrix $[T]_{\widehat{\mathcal{L}},\widehat{\mathcal{L}}} = (I - G_{\widehat{\mathcal{L}},\widehat{\mathcal{L}}})^{-1} $ can be identified.
	
	Now we consider a measured vertex $j \in \mathcal{L}_k$, with $\mathcal{N}_j^- \subseteq \widehat{\mathcal{L}}$. Then, there exists a set of excited vertices $\mathcal{R}_j \subseteq \mathcal{R} \cap \widehat{\mathcal{L}}$ such that the following relation holds:
	\begin{equation} \label{eq:k1}
		[T]_{j,\mathcal{R}_j} = G_{j, \mathcal{N}_j^-} [T]_{\mathcal{N}_j^-,\mathcal{R}_j},
	\end{equation} 
	in which $[T]_{\mathcal{N}_j^-,\mathcal{R}_j}$ is a submatrix of $[T]_{\widehat{\mathcal{L}},\widehat{\mathcal{L}}}$, and it is full row rank as vertex $j$ satisfies \eqref{eq:rankcond}. Then, the entries in $G_{j, \mathcal{N}_j^-}$ are uniquely solved from \eqref{eq:k1}, given the fact that both $[T]_{j,\mathcal{R}_j}$ and $ [T]_{\mathcal{N}_j^-,\mathcal{R}_j}$ are identified from the data.

	In the case that vertex $j$ is unmeasured, $j$ has to be excited. Then, we follow a similar reasoning as the analysis for $j \in \mathcal{L}_2$ and thereby obtain \eqref{eq:TCjN}. If all the in-neighbors of $j$ are excited, then identifiability analysis of $G_{j,\mathcal{N}_j^-}$ follows the same argument as $j \in \mathcal{L}_2$. While $\widehat{\mathcal{N}}^-_j \ne \emptyset$, i.e. there exist unexcited but measured in-neighbors of $j$, we cannot acquire all the elements in the transfer matrix $[T]_{\mathcal{C}_j, {\mathcal{N}}^-_j}$ in \eqref{eq:TCjN}. However, we note that $\mathcal{S}_i$ in this case is the union of $\widehat{\mathcal{N}}^-_j$ and $\bigcup_{i\in {\mathcal{N}}^-_j}^{} \mathcal{N}^+_i$ that satisfies \eqref{eq:thm:rank2}. A minimum disconnecting set $\mathcal{D}$ from $(\mathcal{R}_j \cup \mathcal{S}_i) \setminus j$ to $\mathcal{C}_j$ is also a disconnecting set from $\widehat{\mathcal{N}}^-_j$ to $\mathcal{C}_j$, owing to $\widehat{\mathcal{N}}^-_j \subset \mathcal{C}_j$. Lemma~\ref{lem:factorization} is thereby applied to have
	\begin{equation}
		[T]_{\mathcal{C}_j,\widehat{\mathcal{N}}^-_j} = [T]_{\mathcal{C}_j,\mathcal{D}} Q_3 = [T]_{\mathcal{C}_j,\mathcal{R}_j \setminus j} Q_1^\dagger Q_3,
	\end{equation}
	and thus
	$
	[T]_{\mathcal{C}_j, \mathcal{N}_j^-} 
	= 
	\begin{bmatrix}
		[T]_{\mathcal{C}_j, \mathcal{N}_j^-\setminus \widehat{\mathcal{N}}^-_j}  & [T]_{\mathcal{C}_j,\mathcal{R}_j \setminus j} Q_1^\dagger Q_3
	\end{bmatrix},
	$
	where the vertices in $\mathcal{N}_j^-\setminus \widehat{\mathcal{N}}^-_j$ are excited, and $[T]_{\mathcal{C}_j, \mathcal{N}_j^-\setminus \widehat{\mathcal{N}}^-_j}$ is known from data.
	Then, \eqref{eq:TCjN} can be written as 
	\begin{align} 
		&\begin{bmatrix}
			[T]_{\mathcal{C}_j, \mathcal{N}_j^-\setminus \widehat{\mathcal{N}}^-_j}  & 0
		\end{bmatrix} - [C]_{\star \mathcal{N}_j^-} 
		\\
		&= 
		\underbrace{\begin{bmatrix}
				[T]_{\mathcal{C}_j, j} & [T]_{\mathcal{C}_j,\mathcal{R}_j \setminus j}
		\end{bmatrix}}_{[T]_{\mathcal{C}_j, \mathcal{R}_j}}
		\begin{bmatrix}
			G_{j,\mathcal{N}_j^-} \\ Q_1^\dagger
			\sum_{i \in \mathcal{N}_j^-}  Q_{2i} G_{\mathcal{N}_i^+ \setminus j, i} - Q_1^\dagger Q_3
		\end{bmatrix},
		\nonumber
	\end{align}
	%	
	%	
	%	\begin{equation*}
		%	-\mathbf{e}_i  = \underbrace{\begin{bmatrix}
				%		[T]_{\mathcal{C}_j, j} & [T]_{\mathcal{C}_j,\mathcal{R}_j \setminus j}
				%		\end{bmatrix}}_{[T]_{\mathcal{C}_j, \mathcal{R}_j}}
		%	\begin{bmatrix}
			%	G_{ji} \\ 
			%	Q_1^\dagger Q_2 G_{\mathcal{N}_i^+ \setminus j, i} -  Q_1^\dagger Q_3
			%	\end{bmatrix}.
		%	\end{equation*}
	in which $[T]_{\mathcal{C}_j, \mathcal{R}_j}$ and the two matrices on the left-hand side can be obtained from the measurement data $(y,r)$. Therefore, $G_{j,\mathcal{N}_j^-}$ can be uniquely obtained due to the left invertible transfer matrix $[T]_{\mathcal{C}_j, \mathcal{R}_j}$. 
	
	Consequently, by induction, we have verified that $G_{j,\mathcal{N}_j^-}$ is identifiable, for all $j \in \mathcal{L}_k$ under the conditions of Theorem~\ref{thm:identifiability}.
	This finalizes the proof.
	
	\subsection{Proof of Proposition \ref{pro:thm-necess}}
	\label{ap:pro:thm-necess}
	 {\blue It is implied from \eqref{eq:rem:path1} and \eqref{eq:rem:path2} that 
	 	\begin{align} \label{eq:proofnecess}
	 		b_{\mathcal{R}_j \setminus j \rightarrow \mathcal{C}_j}  = |\mathcal{R}_j| - 1 = b_{(\mathcal{R}_j \cup \mathcal{S}_j) \setminus j \rightarrow \mathcal{C}_j}.
 		\end{align}
	 	Therefore, there exists a set $\mathcal{P}$ of vertex-disjoint paths from $\mathcal{R}_j \setminus j$ to $\mathcal{C}_j \setminus k$ with $|\mathcal{P}| = |\mathcal{R}_j| - 1$ and a measured vertex $k \in \mathcal{C}_j$ that is a descendant of $j$ in the acyclic graph $\mathcal{G}$ such that a path $p_{j \rightarrow k}$ from $j$ to $k$ is vertex-disjoint with $\mathcal{P}$. Note that $\mathcal{P}$ is also a set of vertex disjoint paths from $\mathcal{R}_j \setminus j$ to $\mathcal{C}_j$.
%	 	as well as from $((\mathcal{R}_j \cup \mathcal{S}_j) \setminus j$ to $\mathcal{C}_j$ (due to \eqref{eq:rem:path1}). 
	
		The in-neighbors of the unmeasured vertex $j$ can be divided into two categories: the unexcited but measured vertices $\widehat{\mathcal{N}}_j^- \subseteq \mathcal{C}_j$ and excited vertices $\widecheck{\mathcal{N}}_j^- : = {\mathcal{N}}_j^- \setminus \widehat{\mathcal{N}}_j^-$.
%	there are $|\mathcal{R}_j| - 1$ vertex-disjoint paths from $\mathcal{R}_j \setminus j$ to $\mathcal{C}_j$, where $\mathcal{C}_j$ includes all the measured in-neighbors of $j$, i.e. $\widehat{\mathcal{N}}_j^- \subseteq \mathcal{C}_j$. Let $\mathcal{P}$ be a set of the ,  
	
 First, we prove that there are $|\widehat{\mathcal{N}}_j^-|$ vertex-disjoint paths in $\mathcal{P}$ from $\mathcal{R}_j\setminus j$ to $\widehat{\mathcal{N}}_j^-$. The proof is by contradiction. Assume that there are no $|\widehat{\mathcal{N}}_j^-|$ vertex-disjoint paths from $\mathcal{R}_j \setminus j$ to $\widehat{\mathcal{N}}_j^-$, i.e. $b_{\mathcal{R}_j \setminus j \rightarrow \widehat{\mathcal{N}}_j^- } < |\widehat{\mathcal{N}}_j^- |$.  
 Note that $\mathcal{R}_j \setminus j \subset (\mathcal{R}_j \cup \mathcal{S}_j) \setminus j $, and $b_{(\mathcal{R}_j \cup \mathcal{S}_j) \setminus j \rightarrow \widehat{\mathcal{N}}_j^-} = |\widehat{\mathcal{N}}_j^- |$. Therefore,  we obtain
	 \begin{align*}
	 	 	b_{\mathcal{R}_j \setminus j \rightarrow \mathcal{C}_j} < b_{(\mathcal{R}_j \cup \mathcal{S}_j) \setminus j \rightarrow \mathcal{C}_j},
	 \end{align*}
 	since $b_{(\mathcal{R}_j \cup \mathcal{S}_j) \setminus j \rightarrow \mathcal{C}_j}$ counts not only all the vertex disjoint paths from  $\mathcal{R}_j \setminus j$ to $\mathcal{C}_j$ but also $|\widehat{\mathcal{N}}_j^-|$ vertex disjoint paths from $\widehat{\mathcal{N}}_j^- \subseteq \mathcal{S}_j \setminus j$ to itself in $\mathcal{C}_j$. 
  This then leads to a contradiction to the condition \eqref{eq:proofnecess}. 
 	
 	Next, denote $\mathcal{R}_N : = (\mathcal{R}_j \setminus j) \cap \widecheck{\mathcal{N}}_j^-$. If $\mathcal{R}_N = \varnothing$, then equation  \eqref{eq:pro:cond1} obviously holds true due to $b_{\mathcal{R}_j \setminus j \rightarrow \widehat{\mathcal{N}}_j^- } = |\widehat{\mathcal{N}}_j^- |$. Let's discuss the case where $\mathcal{R}_N \ne \varnothing$. In this case, we prove that there are $|\widehat{\mathcal{N}}_j^-|$ vertex disjoint paths in $\mathcal{P}$ from $(\mathcal{R}_j \setminus j) \setminus \mathcal{R}_N$ to $\widehat{\mathcal{N}}_j^-$.  
 	We proceed with the proof by contradiction. 
 	Assuming that $b_{(\mathcal{R}_j \setminus j) \setminus \mathcal{R}_N \rightarrow \widehat{\mathcal{N}}_j^- } < |\widehat{\mathcal{N}}_j^- |$, then, due to $b_{\mathcal{R}_j \setminus j \rightarrow \widehat{\mathcal{N}}_j^- } = |\widehat{\mathcal{N}}_j^- |$ proved above, there must be a path $p_{i \rightarrow \mu} \in \mathcal{P}$ from an excited vertex $i \in\mathcal{R}_N \subseteq \mathcal{R}_j \setminus j$ to a measured vertex $\mu \in  \widehat{\mathcal{N}}_j^-$. On the other hand, we can construct a set $\mathcal{P}_S$ of vertex disjoint paths from $(\mathcal{R}_j \cup \mathcal{S}_j) \setminus j$ to $\mathcal{C}_j$ as follows: 
 	\begin{align*}
 		\mathcal{P}_S =  \mathcal{P} \setminus p_{i \rightarrow \mu}  \cup p_{i \rightarrow k} \cup p_{\mu \rightarrow \mu},
 	\end{align*}
 	where $p_{\mu \rightarrow \mu}$ is the vertex $\mu$ itself, since $\mu$ is included in both $\mathcal{S}_j$ and $ \mathcal{C}_j$, and $p_{i \rightarrow k}$ is a path obtained by concatenating the edge $(i,j)$ and the path $p_{j \rightarrow k}$. It is not hard to see that $p_{i \rightarrow k}$, $p_{\mu \rightarrow \mu}$, and the paths in $\mathcal{P} \setminus p_{i \rightarrow \mu}$ are vertex disjoint. Therefore, we obtain 
 	\begin{align*}
 		b_{(\mathcal{R}_j \cup \mathcal{S}_j) \setminus j \rightarrow \mathcal{C}_j} = |\mathcal{P}_S| = |\mathcal{P}| + 1 > b_{\mathcal{R}_j \setminus j \rightarrow \mathcal{C}_j},
 	\end{align*}
 	which leads to a contradiction to \eqref{eq:proofnecess}. 
 	
 	Finally, since all the vertices in $\widecheck{\mathcal{N}}_j^-$ are excited, and $b_{(\mathcal{R}_j \setminus j) \setminus \mathcal{R}_N \rightarrow \widehat{\mathcal{N}}_j^- } = |\widehat{\mathcal{N}}_j^- |$ with $(\mathcal{R}_j \setminus j) \setminus \mathcal{R}_N \cap \widecheck{\mathcal{N}}_j^- = \varnothing$, the equation \eqref{eq:pro:cond1} holds.
%  
%
%
% 	 
% 	
% 	Since , the path from $i$ to $k \in \mathcal{C}_j$ via $j$ is vertex disjoint with $\mathcal{P}$.
% 	
% 		From \eqref{eq:rem:path1}, there exists a measured vertex $k \in \mathcal{C}_j$ that is a descendant of $j$ in the acyclic graph $\mathcal{G}$ such that a path from $j$ to $k$ is vertex-disjoint with $\mathcal{P}$.
	}

	\subsection{Proof of Theorem \ref{thm:treecover}}
	\label{ap:thm:treecover}
	
	{\blue Before we prove the theorem, we provide the following lemma that we will use.
	\begin{lemma}
		\label{lem:thmtreecover}
		Consider a network model set $\mathcal{M}$ associated with an acyclic graph $\mathcal{G}$ with the vertex set $\mathcal{V} = \mathcal{R} \cup \mathcal{C}$. If condition (a) in Theorem~\ref{thm:treecover} is satisfied, then the following statements hold.
		\begin{enumerate}
			\item For any measured vertex $j \in \mathcal{V}$, all the edges incident from $\mathcal{N}_j^-$ to $j$ belong to distinct trees. 
			
			\item A vertex $j \in \mathcal{V}$ can have more than one in-neighbor only
			if it is measured.
			
			\item If a vertex $j \in \mathcal{V}$ is excited and not measured, then $|\mathcal{N}_j^-| \leq 1$.
		\end{enumerate}
	\end{lemma}
	\begin{proof}
		The first statement is proved by contradiction. Suppose that two edges incident from $\mathcal{N}_j^-$ to $j$ are included within the same tree $\mathcal{T}$. In such a case, we can identify two distinct paths from the root of $\mathcal{T}$ to $j$ in $\mathcal{T}$. This contradicts the definition of trees, where only one directed path is allowed from the root of a tree to every other vertex in the tree. Thus, the edges incident from $\mathcal{N}_j^-$ to $j$ should be included in distinct trees.
	 	
	 	From the definition of tree covering, all edges incident from $j$ to its out-neighbors are included in the same tree. Therefore, if $j$ has multiple in-neighbors, then $j$ must be an internal vertex of a tree as well as a leaf of another tree in the tree covering. The condition (a) in Theorem~\ref{thm:treecover}   requires all the leaves to be measured, and hence  $j$ should be measured.
	 	
	 	Note that the third statement is the contrapositive of the second one and hence holds.
	\end{proof}
	
	Now we proceed to show the statements of Theorem~\ref{thm:treecover},} and first,
	we prove that $\mathcal{M}$ is generically identifiable if condition (a) holds. 
	Let $\mathcal{G}$ be a composition of disjoint trees $\mathcal{T}_1, \mathcal{T}_2, \cdots, \mathcal{T}_\kappa$, where all the roots are excited, all the leaves are measured, and the rest of vertices are either excited or measured. Therefore, all the sources in $\mathcal{G}$ are excited and all the sinks in $\mathcal{G}$ are measured. Furthermore, 
	{\blue it follows from Lemma~\ref{lem:thmtreecover} that} a vertex has more than one in-neighbor \textit{only if} it is measured.

	To prove generic identifiability of $\mathcal{M}$, {\blue we show that the conditions of this theorem imply the graph-based conditions in Corollary~\ref{coro:paths}.} More specifically, we show that \eqref{eq:rem:path0} is satisfied for every measured vertex in $\mathcal{G}$, while for each excited but unmeasured vertex, there exists a measured vertex set $\mathcal{C}_j$ and excited vertex set $\mathcal{R}_j$ including $j$ such that \eqref{eq:rem:path1} and \eqref{eq:rem:path2} hold. 
	
	First, we consider all the measured vertices in the acyclic graph $\mathcal{G}$. 
	It is implied in Definition~\ref{def:disjtree} that each vertex and its out-neighbors are included in the same tree. Thus, for any vertex $j \in \mathcal{C}$, all the edges incident from the vertices in $\mathcal{N}_j^-$ to $j$ belong to distinct trees.
	Therefore, there exist at least $|\mathcal{N}_j^-|$ vertex-disjoint paths from the roots of the disjoint trees $\mathcal{T}_1, \mathcal{T}_2, \cdots, \mathcal{T}_\kappa$ to $\mathcal{N}_j^-$.  
	As in condition (a) of the theorem, the root of each tree is excited, we have $b_{\mathcal{R}_\mathcal{T} \rightarrow \mathcal{N}_j^-} = |\mathcal{N}_j^-|$, where $\mathcal{R}_\mathcal{T}$ is the set of all the roots in the trees. As a result, \eqref{eq:rem:path0} holds because of $\mathcal{R}_\mathcal{T} \subseteq \mathcal{R}$.

	{\blue Next, all the unmeasured vertices in $\mathcal{G}$ are analyzed, and we will show that they satisfy the two equations \eqref{eq:rem:path1} and \eqref{eq:rem:path2}. Let $j$ be any unmeasured vertex in $\mathcal{G}$, which is excited due to $\mathcal{R} \cup \mathcal{C} = \mathcal{V}$. If condition (a) of Theorem~\ref{thm:treecover} holds, then according to Lemma~\ref{lem:thmtreecover}, $|\mathcal{N}_j^-| \leq 1$, i.e., $j$ is a source or has only one in-neighbor in $\mathcal{G}$. When $j$ is a source, the two equations \eqref{eq:rem:path1} and \eqref{eq:rem:path2} are obvious due to $\mathcal{S}_j = \varnothing$. When $j$ has single in-neighbor $i$ in $\mathcal{G}$, then the set $\mathcal{S}_j$ in \eqref{eq:rem:path2} is specified as
		\begin{equation}
			\mathcal{S}_j =
			\begin{cases}
				i \cup {\mathcal{N}}_i^+ & \text{if $i$ is unexcited but measured;}\\
				{\mathcal{N}}_i^+ & \text{otherwise.}
			\end{cases}       
	\end{equation}}%
	Furthermore, both vertices $i$ and $j$ are in the same tree, denoted by $\mathcal{T}_0$, and we can find the first measured descendant of $j$ in $\mathcal{T}_0$, denoted by $k$, such that either $(j,k)$ is an edge in $\mathcal{T}_0$, or there are no measured vertices on the directed path from $j$ to $k$. 
	%	The existence of $k$ is guaranteed because each excited vertex has only one in-neighbor if condition (a) holds.
	Hereafter, two cases are treated. 
	
	{\blue 
		\textit{Case I: The in-neighbor $i$ is excited, and there is a unique directed path from $i$ to $k$ via $j$.} In this case, $\widehat{\mathcal{N}}_j^- = \varnothing$ since $j$ does not have any unexcited (and measured) in-neighbors. Let $\mathcal{R}_j = \{j\}$ and $\mathcal{C}_j = \{k\}$, which yields $b_{\mathcal{R}_j \rightarrow \mathcal{C}_j} = 1$, and $b_{\mathcal{S}_j\setminus j \rightarrow \mathcal{C}_j} = b_{{\mathcal{N}}_i^+\setminus j \rightarrow \mathcal{C}_j} = 0$. Therefore, \eqref{eq:rem:path1} and \eqref{eq:rem:path2} are satisfied. 
		%	
		% this path is in $\mathcal{T}_a$, namely, it only consists of edges in $\mathcal{T}_a$. As a result, we directly obtain $G_{ji}$ identifiable from Corollary~\ref{coro:singlepath}. 
		%	
		
		\textit{Case II: The in-neighbor $i$ is excited, and there are multiple paths from $i$ to $k$.} 
%		Let $p_0$ be the path from $i$ to $k$ via $j$ and $\mathcal{P}$ be the set of all the paths from $i$ to $k$. Since there is no measured vertex on the path $p_0$ from $j$ to $k$, except $k$, then each vertex on $p_0$ only has exactly one in-neighbor, except $k$. 
		Let $\bar{\mathcal{N}}_k^-$ be the subset of in-neighbors of $k$ such that for each vertex in $\bar{\mathcal{N}}_k^-$, there is a path from $i$ to $k$ via this vertex. Following the definition of disjoint trees, each edge from a vertex in $ \bar{\mathcal{N}}_k^-$ to $k$ should be assigned to a different tree, and therefore, we can find $|\bar{\mathcal{N}}_k^-|$ disjoint trees containing the common vertex $k$, and we denote the set of these trees by 
		$\Upsilon = \{\mathcal{T}_0, \mathcal{T}_1, ..., \mathcal{T}_\ell\}$, where $\ell = |\bar{\mathcal{N}}_k^-|$. Furthermore, in each tree $\mathcal{T}_m$, $m = 1,2, ..., \ell$, there exists a measured vertex, other than $k$. We show this statement by contradiction. If there is a tree $\mathcal{T}_m$ which contains no measured vertex other than $k$. Then $k$ is the only leaf of this tree, and there will not be a path from $i$ to $k$ passing any internal vertex or the roof of $\mathcal{T}_m$. This gives a contradiction, as $\mathcal{T}_m$ contains a vertex in $\bar{\mathcal{N}}_k^-$.
		
		Let $\bar{\mathcal{C}}_j$ collect one measured vertex in each tree in $\Upsilon \setminus \mathcal{T}_0$ such that $|\bar{\mathcal{C}}_j| = \ell - 1$. 
	Furthermore, let  $\bar{\mathcal{R}}_{j}$ collect all the roots of the trees in $\Upsilon \setminus \mathcal{T}_0$, and hence $|\bar{\mathcal{R}}_{j}| =  |\bar{\mathcal{C}}_{j}|$, and the vertices in $\bar{\mathcal{R}}_{j}$ are excited by condition (a) of Theorem~\ref{thm:treecover}. It is clear that each tree in $\Upsilon \setminus \mathcal{T}_0$ contains a unique path in the tree from its root vertex in $\bar{\mathcal{R}}_{j}$ to a vertex in $ \bar{\mathcal{C}}_{j}$, which 
	 implies $b_{\bar{\mathcal{R}}_{j} \rightarrow \bar{\mathcal{C}}_{j}} = |\bar{\mathcal{R}}_{j}|$ because the trees  are disjoint. Furthermore, these paths are disjoint with the path from $j$ to $k$ in $\mathcal{T}_0$. Then, we choose $\mathcal{R}_j = \bar{\mathcal{R}}_{j} \cup j$ and $\mathcal{C}_j = \bar{\mathcal{C}}_{j} \cup k$ to yield 
%	 {\blue Recall that the removal of $\bar{\mathcal{C}}_{j}$ results in no paths from $\mathcal{N}_i^+$ (excluding $j$) to $k$ within $\mathcal{T}_a$.} Consequently, we have
	%	When vertex $i$ is not reachable from $\bar{\mathcal{R}}_{j}$, the set
	%	$\bar{\mathcal{C}}_j$ is a minimum disconnecting set from $\bar{\mathcal{R}}_{j}$ to $ \bar{\mathcal{C}}_{j} \cup k$ and from $\mathcal{N}_i^+ \setminus j$ to $\bar{\mathcal{C}}_{j} \cup k$, simultaneously.
	%	Consequently, we have
	\begin{align*}
		b_{\mathcal{R}_j \rightarrow \mathcal{C}_j} & = b_{\bar{\mathcal{R}}_{j} \rightarrow \bar{\mathcal{C}}_{j}} + b_{j \rightarrow k} = |\bar{\mathcal{R}}_{j}| + 1 = |\mathcal{R}_j|.
	\end{align*}
	Furthermore, we note that removing $\bar{\mathcal{C}}_{j}$ results in no paths from $\mathcal{N}_i^+$ (excluding $j$) to $k$, and hence
	\begin{align*}
		b_{(\mathcal{R}_j \cup \mathcal{N}_i^+) \setminus j \rightarrow \mathcal{C}_j} & = b_{\mathcal{R}_j \setminus j \rightarrow \mathcal{C}_j} =  |\mathcal{R}_j| - 1.
	\end{align*}
	Thereby, both \eqref{eq:rem:path1} and \eqref{eq:rem:path2} are satisfied, and the network model set $\mathcal{M}$ is generically identifiable.}

	Then the proof regarding the condition~(a) has been completed, and the proof for the second condition then directly follows by considering the transpose graph of $\mathcal{G}$. If condition (b) holds for $\mathcal{G}$, then condition (a) is satisfied in the transpose graph $\tilde{\mathcal{G}}$. Therefore, the model set $\tilde{\mathcal{M}}$ associated with $\tilde{\mathcal{G}}$ is generically identifiable. Therefore, it follows from Lemma~\ref{lem:dualty} that the original model set $\mathcal{M}$ is generically identifiable.	
	That completes the proof. 
	
	{\blue 
			\section*{Acknowledgment}
 	This work has benefited from the valuable feedback and constructive criticism from several anonymous reviewers. We are grateful for their time and expertise, and for their contributions to improving the quality and clarity of our paper. 
		}
	
	\bibliographystyle{abbrv}        
	\bibliography{netid}

\end{document}